\documentclass[aos,seceqn,MSNbibl,dvips]{arximspdf}
\usepackage{mathbh}

\doi{10.1214/12-AOS1011} %kopijuoti is PTS
\volume{40}
\issue{3}
\pubyear{2012}
\firstpage{1524}
\lastpage{1549}

\makeatletter

\newcommand{\one}{{\mathbh1}}

\newcommand{\cal}{\mathcal}

\newcommand{\R}{\mathbb R}
\newcommand{\N}{\mathbb N}

\newtheorem{theorem}{Theorem}[section]
\newproclaim{remark}{Remark}[section]

\newcommand{\CK}{{\cal{K}}}
\newcommand{\CM}{{\cal{M}}}
\newcommand{\CY}{{\cal{Y}}}
\newcommand{\CN}{{\cal{N}}}
\renewcommand{\a}{\alpha}
\newcommand{\e}{\varepsilon}
\renewcommand{\b}{\beta}
\newcommand{\g}{\gamma}

\newcommand{\CF}{{\cal{F}}}
\renewcommand{\CK}{{\cal{K}}}

\makeatother

\begin{document}
\begin{frontmatter}

\title{Minimax signal detection in ill-posed inverse~problems}
\runtitle{Detection in ill-posed inverse problems}

\begin{aug}
\author[A]{\fnms{Yuri I.} \snm{Ingster}\thanksref{t1}\ead[label=e1]{yurii\_ingster@mail.ru}},
\author[B]{\fnms{Theofanis} \snm{Sapatinas}\corref{}\ead[label=e2]{fanis@ucy.ac.cy}}
\and
\author[C]{\fnms{Irina A.} \snm{Suslina}\ead[label=e3]{isuslina@mail.ru}}
\runauthor{Y. I. Ingster, T. Sapatinas and I. A. Suslina}
\affiliation{Saint Petersburg Electrotechnical University,
University of Cyprus
and~Saint Petersburg University of Information Technologies,
Mechanics~and~Optics}
\address[A]{Y. I. Ingster\\
Department of Mathematics II\\
Saint Petersburg Electrotechnical University\\
Prof. Popov str., 5\\
197376, St. Petersburg\\
Russia \\
\printead{e1}}
\address[B]{T. Sapatinas\\
Department of Mathematics\\
\quad and Statistics \\
University of Cyprus \\
P.O. Box 20537 \\
Nicosia 1678\\
Cyprus \\
\printead{e2}}
\address[C]{I. A. Suslina\\
Department of Mathematics \\
Saint Petersburg University of Information\\
\quad Technologies, Mechanics and Optics\\
Kronverkskiy pr., 49\\
197101, St. Petersburg\\
Russia\\
\printead{e3}} %adresu isvedimo komanda gale!
\end{aug}

\thankstext{t1}{Supported in part by RFBI Grant 11-01-00577
and by Grant NSh-4472.2010.1.}

% HISTORY:
\received{\smonth{5} \syear{2011}}
\revised{\smonth{4} \syear{2012}}

% ABSTRACT
%
\begin{abstract}
Ill-posed inverse problems arise in various scientific fields. We
consider the signal detection problem for mildly, severely and
extremely ill-posed inverse problems with $l^q$-ellipsoids (bodies), $q
\in(0,2]$, for Sobolev, analytic and generalized analytic classes of
functions under the Gaussian white noise model. We study both rate and
sharp asymptotics for the error probabilities in the minimax setup. By
construction, the derived tests are, often, nonadaptive. Minimax
rate-optimal adaptive tests of rather simple structure are also constructed.
%
%Ill-posed inverse problems arise in various scientific fields. We
%consider the signal detection problem for mildly, severely and
%extremely ill-posed inverse problems with $l^q$-ellipsoids (bodies)
%for Sobolev, analytic and generalized analytic classes of functions
%under the Gaussian white noise model. We restrict our attention to the
%range $q \in(0,2]$, consisting of the ``standard'' case $q=2$ and the
%``sparse'' case $q \in(0,2)$, that has received considerable
%attention in the nonparametric estimation literature over the last
%decade, with the ``sparse'' case $q \in(0,2)$ considered mostly in
%well-posed problems. We study both rate and sharp asymptotics for the
%error probabilities in the minimax setup. By construction, however,
%the derived tests are, often, nonadaptive. In such cases, minimax
%rate-optimal adaptive tests of rather simple structure are also
%constructed.
%
%Possible extensions to a wider range of sets under the alternative
%hypothesis, than the ones considered in this work, are also discussed
%in some detail.
\end{abstract}

% KEYWORDS
%

\end{frontmatter}

%s1 #&#
\section{Introduction}
\label{introFanis}

We consider the detection problem in linear operator equations from
noisy data. More precisely, we
consider the Gaussian white noise model (GWNM)
%
%e1.1 #&#
%
\begin{equation}
\label{101} dY_{\e}(t)=Af(t)\,dt + \e \,dW(t),\qquad t \in D,
\end{equation}
where $A\dvtx {\cal H} \mapsto L^2(D)$ is a known linear bounded operator,
$ {\cal H} \subset L^2(D)$, \mbox{$D \subset\R$}, $W$ is a standard Wiener
process on $D$, $\e>0$ is a small parameter (the noise level) and $f
\in L^2(D)$ is
the unknown response function (that one needs to detect or
estimate).

We consider below the case where $A$ has a kernel with the
\textit{singular value decomposition} (SVD) $A(t,x) = \sum_{k \in\N}
b_k\psi_k(t)\varphi_k(x)$, in the sense of
\[
(Af) (t)=\int_{D} A(t,x) f(x) \,dx=\sum
_{k \in\N} b_k\psi_k(t) \int
_{D} f(x) \varphi_k(x) \,dx,\qquad t \in D,
\]
with $b_k>0$, $k \in\N$, and orthonormal bases $\{\psi_k\}_{k\in\N}$
and $\{\varphi_k\}_{k\in\N}$. (Here,
$\N=\{1,2,\ldots\}$ is the set of natural numbers.)

Thus, the GWNM~(\ref{101}) generates an equivalent discrete
observational model in the
sequence space, called the Gaussian sequence model (GSM),
%
%e1.2 #&#
%
\begin{equation}
\label{103} y_k=b_k\theta_k+\e
\xi_k,\qquad k\in\N,
\end{equation}
where $y_k = \int_{D} \psi_k(t) \,dY_{\e}(t)$, $b_k >0 $, $\theta_k =
\int_{D} f(t) \varphi_k(t) \,dt$, $\e>0$, and $\xi_k \stackrel{
\mathrm{i.i.d.}}{\sim} \CN(0,1)$, $k\in\N$.
%(Note also that, by Parseval's equality, $\theta=\{\theta_k\}_{k\in\N}
The effect of the ill-posedness of the inverse problem is
clearly seen in the decay of $b_k$ as $k \to\infty$. As $k \to
\infty$, $b_k \theta_k$ usually gets weaker and is then more
difficult to detect or estimate $\theta=\{ \theta_k\}_{k\in\N}$.

The GSM~(\ref{103}) can be rewritten in the (equivalent) form
%
%e1.3 #&#
%
\begin{equation}
\label{104} x_k=\theta_k+\e\sigma_k
\xi_k,\qquad k\in\N,
\end{equation}
where $x_k=y_k/b_k$ and $\sigma_k=b_k^{-1}>0$, $k\in\N$. In this
situation, the difficulty of ill-posedness, and hence any asymptotic
results, is measured by the rates (type of growth) of $\sigma_k$ as
$k\to\infty$. For polynomial rates, that is, $\sigma_k\asymp
k^{\beta}$, \mbox{$\beta> 0$}, the inverse problem is
called \textit{mildly} (\textit{or softly}) \textit{ill-posed}, for exponential rates,
that is, $\sigma_k\asymp\exp(\beta k)$, $\beta>0$,
is called \textit{severely ill-posed}, and for the case where
$\sigma_{k+1}/\sigma_k \to\infty$ as $k\to\infty$, is called
\textit{extremely ill-posed}. Note that an extremely ill-posed inverse
problem includes power-exponential rates, that is,
$\sigma_k\asymp\exp(\beta k^{\gamma})$, $\beta>0$, $\gamma
>1$.\setcounter{footnote}{1}\footnote
{The relation $c_{n}
\asymp d_{n}$ means that there exist constants $0 <C_1 \leq C_2 <\infty$
and $n_0$ large enough such that $C_1 \leq c_n/d_n \leq C_2$ for $n
\geq n_0$. We say that $c_{n}(\kappa) \asymp d_{n}(\kappa)$
uniformly over $\kappa\in\CK$, if the similar inequalities hold
true for all $\kappa\in\CK$ with constants $0 <C_1 \leq C_2
<\infty$ and $n_0$ which do not depend on $\kappa$. The relation
$c_{n} \sim d_{n}$ means that for any $\delta\in(0,1)$ there exists
$n_0$ large enough such that $1-\delta\leq c_n/d_n \leq1+\delta$
for $n \geq n_0$. The uniform version of the relation
$c_{n}(\kappa) \sim d_{n}(\kappa), \kappa\in\CK$, is defined
similarly. Similar notation is used when $0 <\e\leq\e_0$ for $\e_0$
small enough.}

An important element of the GSMs
(\ref{103}) and~(\ref{104}) is the prior information about the sequence
$\theta=\{\theta_k\}_{k\in\N}$. Successful detection or estimation
of the sequence $\theta=\{\theta_k\}_{k\in\N}$ is possible only if
its elements $\theta_k$, $k\in\N$, tend to zero sufficiently fast as
$k$ tends to infinity, meaning that
$f$ in the GWNM~(\ref{101}) is sufficiently
smooth. A~standard smoothness assumption on $f$ is to assume
that the sequence $\theta=\{\theta_k\}_{k\in\N}$ belongs to an
$l^q$-ellipsoid (body), $0<q<\infty$, in $l^2$, of semi-axes $L/a_k$,
$k\in\N$, that is,
%
%e1.4 #&#
%
\begin{equation}
\label{fanisellips} \tilde\Theta=\tilde\Theta_q(a,L)= \biggl\{\theta
\in l^2\dvtx \sum_{k\in\N} |a_k
\theta_k|^q\leq L^q \biggr\},
\end{equation}
where $a=\{a_k\}_{k\in\N}$, $a_k > 0$, $a_k\to\infty$ as
$k\to\infty$ and $L>0$. [Note that the requirement $a_k > 0$,
$a_k\to\infty$ as $k\to\infty$, ensures that $\tilde\Theta_q(a,L)$
is a compact subset of $l^2$.] The sequence $a=\{a_k\}_{k\in\N}$
characterizes the ``shape'' of the ellipsoid while the parameter $L$
characterizes its ``size.'' This means that for large values of $k$,
the elements $\theta_k$, $k\in\N$, will decrease in $k$ and, hence,
will be small for large $k$. %(Certainly, the set
%$\tilde\Theta_q(a,L)$ is a ball of the radius $L$ in $l^2$ with
%respect to the norm $|\theta|_{a,q}=\left(\sum_{k\in\N}
%|a_k\theta_k|^q\r)^{1/q}$.)
In what follows, we consider minimax signal detection in ill-posed problems
with $l^q$-ellipsoids for the range $q \in(0,2]$.\vadjust{\goodbreak}

The functional sets of the form~(\ref{fanisellips}) that are often
used in various ill-posed inverse problems are the Sobolev classes
of functions (see~\cite{Wah}) and the classes of analytic functions;
see~\cite{Ib}. We also consider a class of generalized analytic
functions. Then, $\tilde\Theta$ in~(\ref{fanisellips}) takes,
respectively, the form
\begin{eqnarray*}
\mathcal{W}_q(\alpha,L) &=& \biggl\{\theta\in l^2\dvtx \sum
_{k \in\N} k^{\alpha q} |\theta_k|^q
\leq L^q \biggr\},
\\
\mathcal{A}_q(\alpha,L) &=& \biggl\{\theta\in l^2\dvtx \sum
_{k \in\N} e^{\alpha k q} |\theta_k|^q
\leq L^q \biggr\},
\\
\mathcal{G}_q(\tau, \alpha,L) &=& \biggl\{\theta\in l^2\dvtx
\sum_{k
\in\N
} e^{\alpha k^{\tau} q} |\theta_k|^q
\leq L^q \biggr\}
\end{eqnarray*}
for some $\alpha>0$, $\tau\geq1$ (the case $\tau=1$ corresponds to
the class of analytic functions) and $L >0$.
%The Sobolev class of functions has been associated, e.g., with the
%estimation of
%derivatives of smooth functions while the classes of analytic and
%generalized analytic
%functions has been used in, e.g., the estimation of the
%initial or boundary conditions in partial differential equations (see,
%e.g.,~\cite{cav,Cav-Tsyb,Cav-Gol-Lep-Tsyb}).

Despite the growing number of works for the estimation problem in
ill-posed inverse problems under the GWNM~(\ref{101}) (see, e.g.,
\cite{Cav-Tsyb,Cav-Gol-Lep-Tsyb,Gol-Khas2} and~\cite{Gol-Khas1}), very
little work exists for the corresponding detection problem; see
Section~4.3.3 of~\cite{IS02} and~\cite{Er} (although their results are obtained
from models that are neither formulated nor immediately seen as
particular ill-posed inverse problems), and~\cite{ISS11} for a problem
related to the Radon transform (see also, the supplementary material
\cite{ISS12}, Remark~6.1). Our aim is to present a~general
framework for the minimax detection study of the aforementioned
ill-posed inverse problems.
%In particular, both rate and sharp
%asymptotics for the error probabilities in the minimax setup are
%studied in detail and minimax tests are constructed. By
%construction, however, the derived tests are, often, nonadaptive.
%In such case, minimax rate-optimal adaptive tests of rather simple
%structure for the various ill-posed inverse problems under study
%are also constructed.
(Nonasymptotic minimax rates of testing for some of the ill-posed
inverse problems under consideration were recently studied
in~\cite{LLM}.)\looseness=-1

The rest of the paper is organized as follows. The general
statement of minimax signal detection in ill-posed inverse problems
is given in Section~\ref{minimax}, while a short
description of the main results and a comparison with similar
results obtained in the corresponding estimation problems are
presented in Section~\ref{Mainres}. The general methods for the
study of minimax signal detection in ill-posed inverse problems with
$l^2$-ellipsoids are given in Section~\ref{subsecgenresF}. In Sections
\ref{R1}--\ref{Smildlyq}, we provide a complete treatment to the
minimax signal detection problem for mildly, severely and extremely
ill-posed inverse problems with $l^q$-ellipsoids, $q \in(0,2]$, for
Sobolev, analytic and generalized analytic classes of
functions under the GSM~(\ref{103}). We study
both rate and sharp asymptotics for the error probabilities in the
minimax setup. By construction, the derived tests are, often, nonadaptive.
In Section~\ref{Ad}, for the ill-posed inverse problems under
consideration, we also construct minimax rate-optimal adaptive tests
of rather simple structure.
%We pinpoint that sharp, rate and rate-adaptive optimality results for
%the case of mildly ill-posed inverse problems with $l^q$-ellipsoids,
%$q \in(0,2)$, for Sobolev classes of functions are of different
%nature, more complicated, and do not follow directly from the general
%framework considered above. However, as we reveal in Section
%unknown link with the sequence models considered in this paper and
%results obtained in other contexts and presented in Chapters 4, 6 and
%7 of~\cite{IS02}. Since these results are scattered in the cited
%reference, and not immediately seen as ill-posed inverse problems, for
%completeness and an immediate access, we formulate and present them in
%Section~\ref{Smildlyq}.
The proofs along with other relevant material can be found in
the supplementary material~\cite{ISS12}.\footnote{Some numbering that
appears in the text
corresponds to numbering in the supplementary material~\cite{ISS12},
Sections 6--11. Also, some
of the references that appear in the reference list are cited in
the supplementary material~\cite{ISS12}.}

%s2 #&#
\section{Signal detection in the GSM: The minimax framework}\label{minimax}

Consider the GSM~(\ref{103}). In order to avoid having a trivial
minimax hypothesis testing problem (i.e., trivial power), one usually
needs to remove a neighborhood around the functional parameter under
the null hypothesis and to add some additional constraints, that are
typically expressed in the form of some regularity conditions, such as
constraints on the derivatives, of the unknown functional parameter of
interest (see, e.g.,~\cite{IS02}, Sections 1.3--1.4).

In view of the above observation, the main object of our study is the
hypothesis testing problem
%
%e2.1 #&#
%
\begin{equation}
\label{13} %\begin{cases}
H_0\dvtx \theta=0 \quad\mbox{versus}\quad
H_1\dvtx \sum_{k\in\N}|a_k
\theta_k|^q\leq1, \qquad\sum_{k\in\N}
\theta_k^2\ge r_\varepsilon^2, %\end{cases}
\end{equation}
where $\theta=\{\theta_k\}_{k\in\N} \in l^2$, $a=\{a_k\}_{k\in\N}$,
$a_k > 0$, $a_k\to\infty$ as $k\to\infty$, $r_\e>0$, \mbox{$r_\e\to
0$}, is a given family, and $q\in(0,2]$. It means that the
set under the alternative corresponds to an $l^q$-ellipsoid
of semi-axes $1/a_k$, $k\in\N$, with an $l^2$-ball of radius $r_\e$
removed. [For simplicity, in
subsequent sections, we focus attention on ellipsoids of the form
(\ref{fanisellips}) with ``size'' $L=1$.]

Consider the sequence $\eta=\{\eta_k\}_{k\in\N}$,
$\eta_k=b_k \theta_k=\theta_k/\sigma_k$, $k\in\N$. Recall that, in
the ill-posed inverse problems under consideration, $\sigma_k=1/{b_k}
\to\infty$ or $b_k
\to0$, as $k\to\infty$. Hence, $\eta\in l^2$, and the GSM~(\ref{103})
is of
the form %\nref{ffgsm}. % with $\eta$ in place of $s$, i.e.,
%
%e2.2 #&#
%
\begin{equation}
\label{ffgsm} y_k=\eta_k+\e\xi_k,\qquad k\in\N.
\end{equation}
%
%where the sequence $\eta=\{\eta_k\}_{k\in\N} \in l_2$, $\xi_k
Thus,~(\ref{13}) can also be written in the
following equivalent form:
%
%e2.3 #&#
%
\begin{equation}
\label{14} %\begin{cases}
H_0\dvtx \eta=0\quad \mbox{versus}\quad H_1\dvtx
\eta\in\Theta_q(r_\e), %\end{cases}
\end{equation}
where the set under the alternative is
determined by the constraints
%
%e2.4 #&#
%
\begin{eqnarray}
\label{ell} \Theta_q&=&\biggl\{\eta\in l^2\dvtx \sum
_{k\in\N}|a_k\sigma_k\eta_k|^q
\leq1\biggr\},\nonumber\\[-8pt]\\[-8pt]
\Theta_q(r_\e)&=&\biggl\{\eta\in
\Theta_q\dvtx \sum_{k\in\N}\sigma_k^2
\eta_k^2\ge r_\varepsilon^2\biggr\};\nonumber
\end{eqnarray}
%
%which corresponds to the norm
%$|\theta|^2=\sum_{k\in\N}\sigma_k^2\eta_k^2$,
that is, the set under the alternative corresponds to an
$l^q$-ellipsoid of semi-axes $1/(a_k\sigma_k)$, $k\in\N$,
with an $l^2$-ellipsoid of semi-axes
$r_\varepsilon/\sigma_k$, $k\in\N$, removed. %(Note that the

We are therefore interested in the minimax efficiency of the
hypothesis testing problem~(\ref{14}) and~(\ref{ell}) for a given
family of sets $\Theta_\e=\Theta_q(r_{\e})\subset l^2$. It is
characterized by asymptotics, as $\e\to0$, of the minimax error
probabilities in the problem at hand. Namely, for a (randomized)
test $\psi$ (i.e., a measurable function of the observation
$y=\{y_k\}_{k\in\N}$ taking values in $[0,1]$), the null hypothesis
is rejected with probability $\psi(y)$ and is accepted with
probability $1-\psi(y)$. Let $P_{\e,\eta}$ be the probability
measure for the GSM~(\ref{ffgsm}), and denote by
$E_{\e,\eta}$ the expectation over this probability measure. Let
$\a_{\e}(\psi) = E_{\e,0}\psi$ be its type I error probability, and
let $\b_\e(\Theta_{\e},\psi)%=E_{\e,\eta}(1-\psi),
= \sup_{\eta\in
\Theta_{\e}}\b_\e(\eta,\psi) $,
$\b_{\e}(\eta,\psi)=E_{\e,\eta}(1-\psi)$, be its maximal type II error
probability. We consider two criteria of asymptotic optimality:

(1) The first one corresponds to the classical Neyman--Pearson
criterion. For $\a\in(0,1)$, we set
$
\b_{\e}(\Theta_{\e},\a) = \inf_{\psi\dvtx
\a_{\e}(\psi)\le\a}\b_{\e}(\Theta_{\e},\psi).
$
We call a family of tests $\psi_{\e,\a}$ \textit{asymptotically
minimax} if
$
\a_{\e}(\psi_{\e,\a})\le\a+o(1),
\b_{\e}(\Theta_{\e},\psi_{\e,\a}) = \b_{\e}(\Theta_{\e},\a)+o(1),
$
where $o(1)$ is a family tending to zero (here, and in what
follows, unless otherwise stated, all limits are taken as $\e
\rightarrow0$).

(2) The second one corresponds to the total error probabilities.
Let\break
$\g_{\e}(\Theta_{\e},\psi)$ be the sum of the type I and the maximal
type II error probabilities, and let $\g_{\e}(\Theta_{\e})$ be the
minimax total error probability, that is,
$
\g_{\e}(\Theta_{\e}) = \inf_{\psi}\g_{\e}(\Theta_{\e},\psi),
$
where the infimum is taken over all possible tests. We call
a~family of tests $\psi_{\e}$ \textit{asymptotically minimax} if
$
\g_{\e}(\Theta_{\e},\psi_{\e})=\g_{\e}(\Theta_{\e})+o(1).
$
It is known that (see, e.g.,~\cite{IS02}, Chapter~2) that
%
%e2.5 #&#
%
\begin{equation}\quad
\label{bg}\quad \b_{\e}(\Theta_{\e},\a)\in[0,1-\a],\qquad
\g_{\e}(\Theta_{\e})=\inf_{\a\in(0,1)}\bigl(\a+
\b_{\e}(\Theta_{\e
},\a)\bigr) \in[0,1].
\end{equation}

We consider the problems of rate and sharp asymptotics for the error
probabilities in the minimax setup. The rate optimality problem
corresponds to the study of the conditions for which
$\g_{\e}(\Theta_{\e})\to1$ and $\g_{\e}(\Theta_{\e})\to0$ and,
under the conditions of the last relation, to the construction of
\textit{asymptotically minimax consistent} families of tests $\psi_{\e}$,
that is,
such that $\g_{\e}(\Theta_{\e},\psi_{\e})\to0$.

%We are interesting in a set $\Theta_{\e}$ of the form
%$\Theta_\e=\Theta(r_{\e})=\{\eta\in\Theta: |\eta|\ge r_\e\}$,
%where $\Theta\subset l_2$ is a given set, $|\cdot|$ is some norm in
%$l_2$ (not necessarily the standard $l_2$-norm) and $r_\e\to0$ is a
%given positive-valued family.
For the set of the form~(\ref{ell}), we use the notation
$\g_\e(\Theta_q(r_\e))=\g_\e(r_\e)$ and
$\b_\e(\Theta_q(r_{\e}),\a)=\b_\e(r_{\e},\a)$, and we are
interested in
the
%(as in, e.g., \nref{ffffalterset}, i.e.,
%$\barTheta_{\e}(r_{\e})=\barTheta_{\beta,\gamma}(r_{\e})$), these
%conditions correspond to some
minimal decreasing rates for the sequence $r_{\e}$ such that
$\g_\e(r_\e)\to0$. Namely, we say that the positive sequence
$r_{\e}^* \rightarrow0$ is a \textit{separation rate} if
%
%e2.6 #&#
%
\begin{eqnarray}
\label{rates1} \g_{\e}(r_{\e})&\to&1,\qquad \b_{\e}(r_{\e},
\a)\to1-\a\qquad\mbox{for any } \a\in(0,1)\nonumber\\[-8pt]\\[-8pt]
&&\eqntext{\mbox{as } r_{\e}/r_{\e
}^*
\to0,}
\\
%
%e2.7 #&#
\label{rates2} \g_{\e}(r_{\e})&\to&0,\qquad \b_{\e}(r_{\e}, \a)\to0
\qquad\mbox{for any } \a\in(0,1)\nonumber\\[-8pt]\\[-8pt]
&&\eqntext{\mbox{as } r_{\e}/r_{\e }^*
\to\infty.}
\end{eqnarray}
In other words, it means that, for small $\e$, one can detect all
sequences $\eta\in\Theta_q(r_{\e})$ if the ratio $r_{\e}/r_{\e
}^*$ is
large, whereas if this ratio is small, then it is impossible to
distinguish between the null and the alternative hypothesis, with
small minimax total error probability. Hence, the rate optimality
problem corresponds to finding the separation rates $r_{\e}^*$ and
to constructing asymptotically minimax consistent families of tests.

On the other hand, the sharp optimality problem corresponds to the
study of the asymptotics of the quantities
$\b_{\e}(\Theta_{\e},\a), \g_{\e}(\Theta_{\e})$ (up to vanishing
terms) and to the construction of asymptotically minimax families of
tests $\psi_{\e,\a}$ and $\psi_{\e}$, respectively. We shall see
(see Section~\ref{subsecgenresF}) that often the sharp asymptotics
are of Gaussian type, that is,
%
%e2.8 #&#
%
\begin{equation}\quad
\label{G} \quad\b_{\e}(r_{\e},\a)=\Phi\bigl(H^{(\a)}-u_{\e}
\bigr)+o(1), \qquad\g_{\e}(r_{\e})=2\Phi(-u_{\e}/2)+o(1),\vadjust{\goodbreak}
\end{equation}
where $\Phi$ is the standard Gaussian distribution function, and
$H^{(\a)}$ is its $(1-\a)$-quantile, that is,
$\Phi(H^{(\a)})=1-\alpha$.
%For the sets $\Theta_\e=\Theta(r_\e)$,
%we use the notation $u_\e(\Theta_\e)=u_\e(r_\e)$.
The quantity $u_\e=u_\e(r_\e)$ is the value of the specific extreme
problem~(\ref{D1}) on the sequence space $l^2$, and the extreme
sequence of this problem determines the structure of the asymptotically
minimax families of tests $\psi_{\e,\a}$ and $\psi_{\e}$. Moreover,
we shall see that if $u_\e(r_\e)\to\infty$, then $\g_\e(r_\e)\to
0,
\b_\e(r_\e,\a)\to0$, and if $u_\e(r_\e)\to0$, then $\g_\e
(r_\e)\to
1, \b_\e(r_\e,\a)\to1-\a$ for any $\a\in(0,1)$; that is, the family
$u_{\e}(r_{\e})$ characterizes \textit{distinguishability} in the
testing problem. The separation rates~$r_{\e}^*$ are usually
determined by the relation $u_{\e}(r_{\e}^*)\asymp1$ (see, e.g.,
\cite{I93,IS02}). Hence, sharp and rate optimality
problems correspond to the study of the extreme problem~(\ref{D1})
and of the asymptotics of the family $u_{\e}(r_\e)$.

%hypothesis testing problem \nref{1.4}-\nref{ell} is of the form

%Hereafter, the relations $A_{\e} \sim B_{\e}$ and $A_{\e} \asymp B_{
%to those mentioned in Section~\ref{introFanis}.
%Let $\one_{\{A\}}$ be the indicator function of a set $A$ and let
%$(a)_+=\max\{0,a\}$.

%s3 #&#
\section{Minimax signal detection in ill-posed inverse problems: A
short description of some of the main
results}\label{Mainres}

Sharp and rate optimality results for the specific ill-posed inverse
problems under consideration are discussed in detail in Section~\ref{Secspec}.
We give below a short description
of the corresponding results for mildly and severely ill-posed inverse
problems with $l^q$-ellipsoids, $q \in(0,2]$, for Sobolev and analytic
classes of functions.

We consider the hypothesis
testing problem~(\ref{14}) and~(\ref{ell}) in the GSM (\ref{ffgsm}).
For the
``standard'' case
$q=2$, consider the extreme problem
%
%e3.1 #&#
%
\begin{equation}
\label{D1} u_\varepsilon^2=u_\e^2(r_\e)=
\frac{1}{2\e^4}\inf_{\eta\in\Theta
(r_\e
)}\sum_{k\in\N}
\eta_k^4.
\end{equation}
Suppose that $\Theta(r_\e)\not=\varnothing$ and $u_\e>0$, and let
there exist an extreme sequence $\{\tilde{\eta}_k\}_{k\in\N}$ in the
extreme problem~(\ref{D1}).
%(Observe the uniqueness of a nonnegative extreme sequence
%$\{\tilde{\eta}_k\}$, ${k\in\N}$, because, by passing to the sequence $
%with elements $z_k=\tilde{\eta}_k^2$, ${k\in\N}$, we obtain the
%minimization problem of a strictly convex function under linear
%constraints.)
Denote\footnote{The values of $\tilde\eta_k$, $w_k$, $k \in\N$, and
$w_0$ depend on $\e$, that is, $\tilde\eta_k= \tilde\eta_{k,\e}$,
$w_k=w_{k,\e}$, $k \in\N$, and $w_0=w_{0,\e}$.}
%
%e3.2 #&#
%
\begin{equation}
\label{seq1} w_k=\frac{\tilde{\eta}_k^2}{\sqrt{2\sum_{k\in\N}\tilde{\eta
}_k^4}},\qquad k\in\N,\qquad w_0=
\sup_{k\in\N} w_k,
\end{equation}
and consider the families of test statistics and tests
%
%e3.3 #&#
%
\begin{equation}
\label{test1} t_\e=\sum_{k\in\N}w_k
\bigl((y_{k}/\e)^2-1 \bigr),\qquad \psi_{\e,H}=
\one_{\{t_\e>H\}},
\end{equation}
%
%(Note that the values of $\tilde\eta_k$, $w_k$, $k \in\N$, and
%$w_0$ depend on $\e$, i.e., $\tilde\eta_k= \tilde\eta_{k,\e}$,
%$w_k=w_{k,\e}$, $k \in\N$, and $w_0=w_{0,\e}$.)
where $\one_{A}$ denotes indicator function of a set $A$.

%t1 #&#
%
\begin{table}
\caption{The asymptotics $u_\e(r_\e)$ as $r_\e\to0 $}\label{tabl1}
\begin{tabular*}{\tablewidth}{@{\extracolsep{\fill}}lcc@{}}
\hline
& \multicolumn{1}{c}{\textbf{Sobolev classes}} & \multicolumn{1}{c@{}}{\textbf{Analytic classes}} \\
\textbf{{Detection problem}} & \multicolumn{1}{c}{$\bolds{(a_k = k^{\alpha})}$}
& \multicolumn{1}{c@{}}{$\bolds{(a_k = \exp\{\alpha k\})}$} \\
\hline
Mildly ill-posed ($\sigma_k=k^{\beta}$)
&$c_1\e^{-2}r_\e^{(4\alpha+4\beta+1)/2\alpha}$&$ c_2\e^{-2}r_{\e}^2
(\log
r_{\e}^{-1})^{-2\beta-1/2} $\\[2pt]
Severely ill-posed ($\sigma_k = \exp\{\beta k\}$)& $\e^{-2}r_\e
^{2} e^{-2\beta
r_\e^{-1/\a}}$&$
\e^{-2}r_\e^{2(\a+\b)/\a}$\\
\hline
\end{tabular*}
\end{table}

The key tool for the study of the above mentioned hypothesis testing
problem is (the general) Theorem~\ref{testing}. It shows that the
family $u_\e=u_\e(r_\e)$ determines distinguishability in the
problem; if $w_0=o(1)$, then it also determines the sharp
asymptotics~(\ref{G}). The rate optimal tests correspond
to a~weighted $\chi^2$-statistic $t_\e$ of the form~(\ref{test1}).
The main reason is that, for $q=2$,\vadjust{\goodbreak} the problem~(\ref{D1}) is quadratically
convex or can be reduced to a convex problem.\footnote{Rate optimal
tests of simpler structure are also given in the supplementary material~\cite{ISS12},
Sections 7--9.}
The key idea is that the $\chi^2$-distance between $P_0$
the probability measure of the observations under the null hypothesis
and the mixture~$P_{\pi_\eta}$ over the product symmetric two-points prior
$\pi_\eta$ for a point $\eta\in\Theta(r_\e)$ is characterized by the
quantity $\sum_{k\in\N}\eta_k^4$, which leads to the extreme problem
(\ref{D1}). Moreover, if $w_0=o(1)$, then the Bayesian log-likelihood
ratio is asymptotically Gaussian under $P_0$, that is,
\[
\log(dP_{\pi_{\tilde\eta}}/dP_0)=-u_\e^2/2+u_\e
t_\e+\delta_\e,
\]
where $u_\e$ is given by~(\ref{D1}), $t_\e$ is given by (\ref
{test1}) and
is asymptotically standard Gaussian and $\delta_\e\to0$ in
$P_0$-probability. This yields the lower bounds of the Gaussian
form. On the other hand, the choice of the optimal coefficients of the test
statistic $t_\e$ given by~(\ref{test1}) leads to a maximin problem that is
reduced to the extreme problem~(\ref{D1}) by convexity arguments.
This yields the corresponding upper bounds. On the other hand,
if $w_0\not= o(1)$, then the test statistic $t_\e$ given by (\ref
{test1}) is not asymptotically Gaussian.
This is due to that fact that, in this case, $a_k$ and $b_k$,
$k \in\N$,
converge fast enough as $k$ increases, implying that only a small number of observations
is important.
[Often, but not always, using embedding properties, these results can
be extended to the case
$q\in(0,2)$; see Remark
\ref{remsparsecase}.]

The asymptotics of the quality of testing
$u_\e(r_\e)$ as $r_\e\to0 $ is presented in Table~\ref{tabl1}, where
$c_1=c_1(\a,\b)>0, c_2=c_2(\a,\b)>0$ are some constants. We
have the sharp asymptotics of the form~(\ref{G}) for mildly
ill-posed inverse problems with Sobolev and analytic classes of
functions, while\vspace*{1pt} the derived asymptotically minimax tests are
based on weighted $\chi^2$-statistics with weights $w_k\ge0,
\sum_k w_k^2=1/2$ and are determined by the extreme problem~(\ref{D1}).
For severely ill-posed inverse problems, however, we do not
have sharp asymptotics of minimax error probabilities but, instead, we get
distinguishability conditions [i.e., $\g_\e(r_\e)\to0$ if and only if
$u_\e(r_\e)\to\infty$ and $\g_\e(r_\e)\to1$ iff $u_\e(r_\e
)\to0
$]. The main reason is that the weights are not
``uniformly small,'' that is, there exist a few coefficients $w_k$ that
are bounded away from zero, and, hence, in this case, we do not have
asymptotic Gaussianity.\looseness=-1

Furthermore, the separation rates $r^*_\e$ as $\e\to0$ are
presented in Table~\ref{tabl2}. (Similar nonasymptotic minimax rates are
recently given in
%
%t2 #&#
%
\begin{table}
\caption{The separation rates $r^*_\e$ as $\e\to0$}\label{tabl2}
\begin{tabular*}{\tablewidth}{@{\extracolsep{\fill}}lcc@{}}
\hline
&\textbf{Sobolev classes} &\textbf{Analytic classes} \\
\textbf{{Detection problem}} & $\bolds{(a_k = k^{\alpha})}$
& $\bolds{(a_k = \exp\{\alpha k\})}$ \\
\hline
Mildly ill-posed ($\sigma_k=k^{\beta}$) &$\varepsilon^{4\alpha
/(4\alpha+4\beta+1)}$&$
\varepsilon(\log\varepsilon^{-1})^{\beta+1/4}$\\
[2pt]
Severely ill-posed ($\sigma_k = \exp\{\beta k\}$)& $((\log
\varepsilon^{-1})/\b)^{-\alpha}$&$
\varepsilon^{\alpha/(\alpha+\beta)}$\\
\hline
\end{tabular*}
\end{table}
\cite{LLM}.) Note that, despite the fact that we have no sharp
asymptotics for
severely ill-posed inverse problems, we get the sharp separation rates
for severely
ill-posed inverse problems with Sobolev $f$, when
$a_k\sim k^\a$ and $b_k\asymp\exp(-\b k)$, $k \in\N$ [i.e.,
$ \g_\e(r_\e)\to0$ as $\lim r_\e/r^*_\e>1$ and $\g_\e(r_\e)\to1$
as $\lim r_\e/r^*_\e<1$].

We do not present the results for extremely ill-posed inverse problems
with generalized analytic $f$ in this short
description because the asymptotics are more complicated in this
case. Roughly speaking, these asymptotics are determined by a
piecewise linear function $u_\e^{\mathrm{lin}}=u_\e^{\mathrm{lin}}(r_\e)$ and,
principally, seem to be of a new type; see Section~\ref{R5} and
remarks therein for details. (Moreover, in subsequent sections, we
consider this case only for $q=2$.)

Consider now the ``sparse'' case $q\in(0,2)$. Then, the results noted
above still hold true for severely
ill-posed inverse problems with Sobolev $f$ or analytic $f$. For mildly
ill-posed inverse problems with analytic $f$, we also get the same
separation rates $r^*_\e$. However, the situation for mildly ill-posed
inverse problems with Sobolev $f$ is more delicate.
%and we refer to the results that one could exclude from~\cite{IS02}.
More precisely, let $\a>0, \b>0$ and set $\lambda=(\alpha+\beta
)/2-\beta/q$. If
$\lambda>0$, then
%one can get from~\cite{IS02}, Sect. 6.2.1
the sharp asymptotics are of the Gaussian type~(\ref{G}) with
\[
u_\e=c_3\e^{-(2\a+1/q-1/2)/(\a+\b(1-2/q))}r_\e^{(2(\a+\b
)+1/q)/(\a+\b(1-2/q))}
\]
for some constant $c_3=c_3(\a,\b,q) >0$, while the separation rates
$r^*_\e$ are of
the form
\[
r_\varepsilon^*= \varepsilon^{(2\alpha+1/q-1/2))/(2(\alpha+\beta)+1/q)}.
\]
The corresponding rate optimal tests are of more complicated structure
and are based on
a different extreme problem; see Section~\ref{Smildlyq}.
%Similar sharp asymptotics hold true for Besov balls under the norm
%$|\cdot|_{\a,q,t}$ when $t\ge q$; if $t< q$, then we have the same
%separation rates, see~\cite{IS02}, Section 6.2.2.

On the other hand, if $\lambda\le0$, then
%one can get from~\cite{IS02}, Sect. 4.4.2
the sharp asymptotics are of the following degenerate type:
\[
\b_\e(\a)=(1-\a)\Phi(-D_\e)+o(1), \qquad\g_\e=
\Phi(-D_\e)+o(1),
\]
where
$
D_\e=n_\e^{-\b}r_\e/\e-\sqrt{2\log(n_\e)}, n_\e=r_\e^{-1/\a},
$
while the separation rates $r^*_\e$ are of
the form
\[
r^*_\e=\Lambda\e^{\a/(\a+\b)}\bigl(\log\bigl(\e^{-1}
\bigr)\bigr)^{\a/2(\a+\b
)}, \qquad\Lambda=\bigl(2/(\alpha+\b)\bigr)^{\a/2(\a+\b)}.
\]
The corresponding rate optimal tests are based on a simple thresholding rule.

%Similar separation rates hold true for Besov balls under the norm
%$|\cdot|_{\a,q,t}$, see~\cite{IS02}, Section 4.4.3.

%These rates results are translated to the Besov and Sobolev balls in
%the functional space of the smoothness $\tilde\a=\a-1/2+1/q$, see

%Similar results hold true for mildly ill-posed problems with Besov
%classes for $q\in(0,2)$

It seems natural to compare the separation rates
$r_\varepsilon^*$ in the detection problem with the minimax accuracy
$R_\varepsilon$
in the corresponding estimation problem using loss functions that
correspond to the norm
which is used for bounding the alternative away from zero. We compare
below the above mentioned minimax rates of testing with the
corresponding minimax rates of estimation. The minimax estimation
problem for the GWNM~(\ref{101}) [or,
equivalently, for the GSM~(\ref{103})] was
studied very intensively in statistical ill-posed inverse problems;
see, for example,~\cite{cav} (and references therein),
\cite{Cav-Tsyb,Cav-Gol-Lep-Tsyb,Gol-Khas2} and~\cite{Gol-Khas1}. The
main\vspace*{1pt} object of the study is the minimax quadratic risk, defined by
$
R_{\e}^2(\CF)=\inf_{\hat f}\sup_{f\in\CF}E_{\e,f}\|\hat{f}-f\|^2,
$
where the infimum is taken over all possible estimators $\hat f$ of
$f$, based on observations from
the GWNM~(\ref{101}).

For the main types of the ill-posed inverse problems and classes of
functions under consideration, with $q=2$, the estimation rates
$R_{\e} =R_{\e}(\CF)$ as $\e\to0$ are presented in Table~\ref{tabl3};
see, for
example,~\cite{cav}. For mildly ill-posed inverse problems with Sobolev
$f$, $q \in(0,2)$, one has, for $\lambda>0$,
\[
R_{\e}\asymp\e^{(\a-1/2+1/q)/(\a+\b+1/q)},
\]
while, for $\lambda<0$, the estimation rates $R_{\e}$
coincide with the separation rates~$r_{\e}^*$ in the corresponding
detection problem; see, for example,
\cite{IS02}, Section~2.8, and references therein. Observe that the
minimax rates of testing are faster than the
corresponding minimax rates of estimation (as it is common in
nonparametric inference, see, e.g.,~\cite{IS02}, Sections 2.10 and
3.5.1), except for the cases of mildly ill-posed inverse problems with
Sobolev $f$, $q\in(0,2)$
and $\lambda\le0$, and the cases of severely ill-posed inverse
problems with Sobolev~$f$ or analytic~$f$, $q=2$. [To the best of our
knowledge, we are not aware of any
minimax estimation results with the case of severely ill-posed
inverse problems for Sobolev~$f$
or analytic $f$, $q \in(0,2)$.]

%
%t3 #&#
%
\begin{table}
\caption{The estimation rates $R_{\e} =R_{\e}(\CF)$
as $\e\to0$}\label{tabl3}
\begin{tabular*}{\tablewidth}{@{\extracolsep{\fill}}lcc@{}}
\hline
&\textbf{Sobolev classes} &\textbf{Analytic classes} \\
\textbf{{Estimation problem}} & $\bolds{(a_k = k^{\alpha})}$
& $\bolds{(a_k = \exp\{\alpha k\})}$ \\
\hline
Mildly ill-posed ($\sigma_k=k^{\beta}$)
&$\varepsilon^{2\alpha/(2\alpha+2\beta+1)}$&$\varepsilon(
\log\varepsilon^{-1})^{\beta+1/2}$\\
[2pt]
Severely ill-posed ($\sigma_k = \exp\{\beta k\}$)& $(\log
\varepsilon^{-1})^{-\alpha}$&$
\varepsilon^{\alpha/(\alpha+\beta)}$\\
\hline
\end{tabular*}
\end{table}

Returning to the signal detection problem, note that, except for the
mildly ill-posed inverse problems with
Sobolev $f$, $q\in(0,2)$ and $\lambda\le0$, or analytic $f$, the
aforementioned separation rates $r_{\e}^*$
still hold true for the \textit{known} parameters $(\a,\b,q)$ associated with
the classes of functions and the ill-posed inverse problems under
consideration. For these cases, the rate-optimal tests also depend on
the parameters
$(\a,\b,q)$. In practice, the parameters $\a$ and $
q$ associated with the considered functional classes are typically
unknown and, very often, the statistician is not confident about the
value of the parameter $\b$ associated with the sequence $b_k$, $k \in
\N$. For unknown parameters $(\a,\b, q)\in\Sigma\subset
\R^2_+\times(0,2]:=(0,\infty)\times(0,\infty)\times(0,2]$, we
have the
so-called \textit{adaptive} problems: in order to distinguish between
the null hypothesis and the ``combined'' alternative, which
corresponds to a~wide enough compact set $\Sigma\subset\R^2_+\times
(0,2]$, it
does not suffice to just require $u_\e=u_\e(r_\e(\a,\b,q),\a,\b
,q)\to
\infty$ for all $(\a,\b,q)\in\Sigma$; instead, one needs that it
should tend to $\infty$ faster than some family
$u_\e^{\mathrm{ad}}\to\infty$, which is a ``payment'' for adaptation; see
\cite{Spok}.

Adaptive rate optimality results for the specific ill-posed inverse
problems under consideration are discussed in
detail in Section~\ref{Ad}. Below, we give a~short description
of these results for mildly and severely ill-posed inverse problems
with Sobolev $f$ or analytic $f$, $q \in(0,2]$.

%
%t4 #&#
%
\begin{table}
\caption{The adaptive separation rates $r_\e^{\mathrm{ad}}$ as $\e\to0$}\label{tabl4}
\begin{tabular*}{\tablewidth}{@{\extracolsep{\fill}}lcc@{}}
\hline
&\textbf{Sobolev classes}&\textbf{Analytic classes}\\
\textbf{{Detection problem}} & $\bolds{(a_k = k^{\alpha})}$
& $\bolds{(a_k = \exp\{\alpha k\})}$ \\
\hline
Mildly ill-posed ($\sigma_k=k^{\beta}$)
&$(\tilde\varepsilon_1)^{4\alpha/(4\alpha+4\beta+1)}$&$
\varepsilon(\log\varepsilon^{-1})^{\beta+1/4}$\\
[2pt]
Severely ill-posed ($\sigma_k = \exp\{\beta k\}$)& $((\log
\varepsilon^{-1})/\b)^{-\alpha}$&
$(\tilde\varepsilon_2)^{\alpha/(\alpha+\beta)}$\\
\hline
\end{tabular*}
\end{table}

For mildly ill-posed inverse problems with Sobolev $f$, $\lambda\le0$,
or analytic $f$, one has $u_\e^{\mathrm{ad}}\asymp1$, while for Sobolev $f$,
with $\lambda> 0$, one has $u_\e^{\mathrm{ad}}=\sqrt{\log\log\e^{-1}}$.
On the other hand, for severely ill-posed inverse problems with
Sobolev $f$ or analytic $f$ (as well as for extremely ill-posed inverse
problems for
generalized analytic~$f$, $q=2$),
one has
$u_\e^{\mathrm{ad}}={\log\log\e^{-1}}$.
These yield the \textit{adaptive separation rates} $r_\e^{\mathrm{ad}}$ as $\e\to0$
presented in Table~\ref{tabl4} (here $q=2$ for mildly ill-posed
problems with
Sobolev $f$), where $ \tilde\e_1=\e\sqrt[4]{\log\log\e^{-1}}$ and
$\tilde\e_2=\e\sqrt{\log\log\e^{-1}}$.
On the other\vspace*{2pt} hand, for the ``sparse'' case $q\in(0,2)$ and $\lambda>0$,
%we have the same adaptive separation rates $r_\e^{\mathrm{ad}}$
%for severely ill-posed problems with $l^q$-ellipsoids for Sobolev and
%analytic classes of functions, while
for mildly ill-posed inverse problems with Sobolev $f$, the adaptive
separation rates $r_\e^{\mathrm{ad}}$ are of the form
%one could exclude from~\cite{IS02}, Section 7.1.3,
%
\[
r_\e^{\mathrm{ad}}= \tilde\e_2^{(2\alpha+1/q-1/2))/(2(\alpha+\beta)+1/q)}.
\]
%
%We have the same adaptive rates for mildly ill-posed problems with
%Besov classes for $q\in(0,2), \lambda>0$.
%
As we shall see
in Section~\ref{Ad}, the rate-optimal adaptive tests are of rather
simple structure for all problems under consideration [except for
the mildly ill-posed problems with Sobolev $f$, $q \in(0,2)$]: they
are based on combinations of tests based on a grid of centered and
normalized statistics of $\chi^2$-type and on simple thresholding.
For mildly ill-posed problems with Sobolev $f$, $q \in(0,2)$, the
rate-optimal adaptive tests are, however, more complicated; see \cite
{IS02}, Chapter 7.
%(see Remark~\ref{remBeqSobFN} in Section~\ref{Ad1q}).

Finally, we mention that, in most cases, the arguments for $q=2$ do not
work for $q \in(0,2)$
[e.g., for mildly ill-posed inverse problems with Sobolev~$f$,
$q\in(0,2)$]; instead,
we have different extreme problems in the space of sequences of
probability measures.
%Moreover, the size and power of the resulting tests are not obtained
%by applying the central limit theorem, since
%there are a few large (bounded away from zero) weights $w_k$, $k\in
%and note that $\sum_{k \in\N} w_k^2=1/2$).
The key ideas for the proofs of the main results for \mbox{$0 < q<2$}, as well
as possible extensions to $0<q\le\infty$ and to a~wider range of
sets under the alternative (i.e., to replace the constraint in the
$\ell^2$-norm by a constraint in the $\ell^p$-norms, $0 < p \leq
\infty$), are discussed in some detail in the supplementary material
\cite{ISS12}, Section 10.

%one can extract the structure of rate-optimal adaptive tests from

%s4 #&#
\section{Minimax signal detection in ill-posed inverse problems: Rate
and sharp asymptotics}
\label{Secspec}

In this section, we consider the GSM~(\ref{ffgsm}) and the hypothesis
testing problem~(\ref{14}) and~(\ref{ell}).

%s4.1 #&#
\subsection{A general result for $l^q$-ellipsoids: The ``standard''
case $q=2$}
\label{subsecgenresF}
%Following the previous discussion, consider
%Consider the Gaussian sequence model \nref{ffgsm}. We are interested
%in the hypothesis testing problem \nref{1.4} with the set under the
%alternative $\Theta_\e=\Theta_2(r_\e)$ given by \nref{ell} with
%$q=2$. Consider now the extreme problem
%u_\varepsilon^2=u_\e^2(r_\e)=\frac{1}{2\e^4}\inf_{\eta\in\Theta(r_\e)}
%Suppose that $\Theta(r_\e)\not=\varnothing$ and $u_\e>0$, and let
%there exist an extreme sequence $\{\tilde{\eta}_k\}_{k\in\N}$ in the
%extreme problem \nref{D1}. (Observe the uniqueness of a nonnegative
%extreme sequence
%$\{\tilde{\eta}_k\}$, ${k\in\N}$, because, by passing to the sequence $
%obtain the minimization problem of a strictly convex function under
%linear constraints.) Denote
%w_k=\frac{\tilde{\eta}_k^2}{\sqrt{2\sum_{k\in\N}\tilde{\eta}_k^4}},
%k\in\N, w_0=\sup_{k\in\N} w_k,
%and consider the following families of test statistics and tests
%t_\e=\sum_{k\in\N}w_k\left((y_{k}/\e)^2-1\r),
%(Note that the values of $\tilde\eta_k$, $w_k$, $k \in\N$, and
%$w_0$ depend on $\e$, i.e., $\tilde\eta_k= \tilde\eta_{k,\e}$,
%$w_k=w_{k,\e}$, $k \in\N$, and $w_0=w_{0,\e}$.)

%The key tool for the study of the hypothesis testing
%problem mentioned in Section~\ref{Mainres} is the following general
%theorem.
%
%th4.1 #&#
%
\begin{theorem}\label{testing}
Let $q=2$, let $u_{\e}$ be determined by the extreme
problem~(\ref{D1}), let the coefficients $w_k$, $k\in\N$, and $w_0$
be as in~(\ref{seq1}) and consider the family of tests $\psi_{\e,H}$
given by~(\ref{test1}).\vspace*{8pt}

\textup{(1) (a)} If $u_\e\to0$, then $\b_\e(r_\e,\a)\to1-\a$ for any $\a
\in
(0,1)$ and $\g_\e(r_\e)\to1$; that is, minimax testing is impossible.
If $u_\e=O(1)$, then $\liminf\b_\e(r_\e,\a)>0$ for any $\a\in(0,1)$
and $\liminf\g_\e(r_\e)>0$; that is, minimax consistent testing is
impossible.

\mbox{}\hphantom{\textup{(1)} }\textup{(b)} If $u_\e\asymp1$ and $w_0=o(1)$, then the family of tests
$\psi_{\e, H}$ of the form~(\ref{test1}) with $H=H^{(\a)}$ and
$H=u_\e/2$ are asymptotically minimax, that is,
\begin{eqnarray*}
\a_{\e}(\psi_{\e,H^{(\a)}})&\leq&\a+o(1),\qquad \b_{\e}\bigl(
\Theta(r_{\e}), \psi_{\e,H^{(\a)}}\bigr)=\b_{\e}(r_{\e},
\a)+o(1),
\\
\g_{\e}\bigl(\Theta(r_{\e}), \psi_{\e,u_\e/2}\bigr)&=&
\g_{\e}(r_{\e})+o(1),
\end{eqnarray*}
%
%
% \g_{\e}(\Theta(r_{\e}), \psi_{\e,u_\e/2})&=&\g_{\e}(r_{\e})+o(1),
and the sharp asymptotics~(\ref{G}) hold true. %, i.e.,
%$$
%$$

(2) If $u_\e\to\infty$, then the family of tests $\psi_{\e,H}$ of
the form~(\ref{test1}) with $H\sim cu_\e$ are asymptotically minimax
consistent for any $c\in(0,1)$, that is,
$\g_\e(\Theta(r_{\e}),\psi_{\e,T_\e})\to0$.
\end{theorem}

The proof is given in the supplementary material~\cite{ISS12}, Section 11.1.

Theorem~\ref{testing} shows that the asymptotics of the quality of
testing is determined by the asymptotics of values $u_\e$ of the the
extreme problem~(\ref{D1}). This latter problem is studied by
using Lagrange multipliers. Then the extreme sequence in the above
mentioned extreme problem is of the form
%
%e4.1 #&#
%
\begin{equation}
\label{seq} \tilde{\eta}_k^2=z_0^2
\sigma_k^2\bigl(1-Aa_k^2\bigr)_+,\qquad
k\in\N,
\end{equation}
where $(a)_+=\max\{0,a\}$, and the quantities $z_0=z_{0,\e}$ and
$A=A_{\e}$ are determined by
the equations
%
%e4.2 #&#
%
\begin{equation}
\label{eq0} \sum_{k\in\N} \sigma_k^2
\tilde{\eta}_k^2=r_\e^2,\qquad \sum
_{k\in\N} a_k^2
\sigma_k^2 \tilde{\eta}_k^2=1.
\end{equation}
%
%&\sum_{k\in\N} a_k^2\sigma_k^2 \tilde{\eta}_k^2=1.
Note that the quantity $A$ determines the \textit{efficient dimension}
$m$ in the specific ill-posed inverse problems considered below: if
$a_k$ is an increasing sequence (it is assumed further), the
efficient dimension is a quantity $m=m_\e\in[1,\infty)$ such that
$Aa_{[m]}^2\le1<Aa_{[m]+1}^2$.

%re4.1 #&#
%
\begin{remark}
\label{asssigmaFN}
Since, in the ill-posed problems under consideration, $\sigma_k\to
\infty$ as
$k\to\infty$, it is immediate that $\sum_{k\in\N}\sigma_k^4=\infty$.
Under this condition, and the fact that
$a_k\to\infty$ as $k\to\infty$,
%a_k\to\infty\mbox{as} k\to\infty,
%If assumption \nref{sigma} is violated, then $J_1$, $J_2$, $J_0$ and
%$z_0$ are
%bounded and from \nref{eq1} we get that that %$u_{\e} \asymp\e^2$,
%leading to the classical
%asymptotic rates.
one can see that, for $r_\e$ small
enough, the equations in~(\ref{eq0}) have a unique solution; see
Proposition 11.2 in the supplementary material~\cite{ISS12}, Section 11.9.
\end{remark}

%re4.2 #&#
%
\begin{remark}
\label{rescalargFF}
Let $u_{\e}=u_{\e}(r_{\e})$ be the value of
the extreme problem~(\ref{D1})
with $a=\{a_k\}_{k\in\N}$ and
$\sigma=\{\sigma_k\}_{k\in\N}$ associated with $\Theta_\e=\Theta
(r_\e)$
given by~(\ref{ell}), and let
$\tilde u_{\e}=\tilde u_{\e}(r_{\e})$
be the corresponding value of the extreme problem similar to (\ref
{D1}) with
$\tilde a=C a=\{C a_k\}_{k\in\N}$ and
$\tilde\sigma=D \sigma=\{D\sigma_k\}_{k\in\N}$ in~(\ref{ell}), for
some positive constants $C$ and $D$.
Then it is easily seen that the relation
$ \tilde u_{\e}(r_{\e}) = (CD)^{-2}u_{\e}(C r_{\e}) $ holds true.
\end{remark}

%s4.2 #&#
\subsection{Application to mildly ill-posed inverse problems with the Sobolev
class of functions}\label{R1}

Consider first the ``standard'' case $q =2$.

%th4.2 #&#
%
\begin{theorem}\label{T1}
%Consider the GSM \nref{ffgsm} and the hypothesis
%testing problem \nref{1.4}-\nref{ell} with $q=2$.
Let $q=2, a_k = k^\alpha$ and $\sigma_k = k^\beta$, $k\in\N$,
$\alpha>0, \beta>0$.

\begin{longlist}[(a)]
\item[(a)] The sharp asymptotics~(\ref{G}) hold with the value $u_\e$ of
the extreme problem~(\ref{D1}) determined by (11.10)
in~\cite{ISS12}.

\item[(b)] The asymptotically minimax family of tests $\psi_{\e,H}$ are
determined by the family of test statistics $t_{\e}$ given by
(\ref{test1}) with $w_k$, $k\in\N$, and $w_0$ as in
(\ref{seq1}), and with $\{\tilde{\eta}_k\}_{k\in\N}$
given by (11.9) with $m$ determined by (11.10) in~\cite{ISS12}.

\item[(c)] The separation rates are of the form
%
%e4.3 #&#
%
\begin{equation}
\label{seprateSobSovF} r_\varepsilon^* = \varepsilon^{4\alpha/(4\alpha
+4\beta+1)}.
\end{equation}
\end{longlist}
\end{theorem}

The proof is given in the supplementary material~\cite{ISS12}, Section 11.3.
%{\rm Under the assumption \nref{sigma}, it is immediately seen that
%the asymptotic results
% obtained in Theorem~\ref{T1} hold true for any $\beta> -1/4$. Hence,
%the separation
% rates \nref{seprateSobSovF} are true not only for the ill-posed
%problems ($\beta>0$)
% under consideration but also for a range of corresponding well-posed
% problems ($\beta\in(-1/4, 0]$). }

%{\rm
%%It is also easy to see that the asymptotic results in Theorem
%%uniformly over
%%$(\a,\b)\in\Sigma$ for any compact set $\Sigma\subset(0, \infty)
%%(-1/4,\infty)$.
%
%%%%%New Revised Remark %%%%
%It is easily seen that the asymptotic results in Theorem~\ref{T1} (as
%well as the asymptotic results in all subsequent theorems of Section
%3) hold true uniformly over
%$(\a,\b)\in\Sigma$ for any compact set $\Sigma\subset\R^2_+ $.
%}

%or ill-posed)
%inverse problems are true.
%
%(a) If $\alpha\leq0$ and $\beta\ge-1/4$, then
%$\gamma_{\e}(\Theta_{\e}) = 1$, i.e., minimax testing is impossible.
%
%(b) If
%$\alpha>0$ and $\beta=-1/4$, then the separation rates are of the form
%$
%r_\varepsilon^* = \varepsilon\log^{1/4}(\e^{-1}).
%$
%
%(c) If $\beta<-1/4$, then we arrive at the classical separation
%rates, i.e., $r_\e^*=\varepsilon$.
%
%4.3.3 in
%3.4.3 in~\cite{IS02}.}
%

%re4.3 #&#
%
\begin{remark}
\label{RemequivSS}
It follows from the evaluations of the functions $J_0$, $J_1$ and~$J_2$
used to express~(\ref{eq0}) (see (11.2) in~\cite{ISS12}) that their
asymptotics are determined by the tails of the sequences $a_k =
k^\alpha
$ and
$\sigma_k = k^\beta$, $k\in\N$, $\alpha>0$,
$\beta>0$. For this reason, in view of Remark~\ref{rescalargFF}, we get
the sharp asymptotics~(11.11) in~\cite{ISS12} for the sequences $a_k \sim k^{\alpha}$ and
$\sigma\sim
k^{\beta}$, $k\in\N$,
$\alpha>0,
\beta>0$, and similar\vadjust{\goodbreak} rate asymptotics for the sequences $a_k \asymp
k^{\alpha}$,
$\sigma\asymp k^{\beta}$, $k\in\N$, $\alpha>0,
\beta>0$. In both cases, the separation rates are still of the form~(\ref{seprateSobSovF}).
%The Remarks~\ref{FFF1},~\ref{FFF1FF} and~\ref{Rem1} still apply to
%these cases too.
%Remark~\ref{FFF1FF} still applies to these cases too.
\end{remark}

%of simple structure}
%
%If $u_\e\to\infty$, then one can construct a family of asymptotically
%minimax consistent
%tests
%of simpler structure than \nref{test1}. Indeed, observe that, by
%r_\e a_{\tilde m}\sim c_1^\alpha, c_1=c_1(\alpha,\beta)> 1,
%u_\e\asymp\frac{r_\e^2}{\e^2\sqrt{\tilde m}\sigma_{\tilde
%m}^2}, {\tilde m}=[m]\in\N
%where $[m]$ is the integral part of $m$. Hence, for an
%integer-valued family $\tilde{m}=\tilde m_\e\to\infty$, one has
%a_{\tilde m+1} r_\e\ge B+o(1), B>1, u_\e\asymp
%%\frac{(r_\e^*)^2}{\e^2\sqrt{ \tilde m}\sigma_{\tilde m}^2}\asymp1,
%For each $m\in\N$, consider the following families of test statistics
%and tests
%t_{\e, m}=\frac{1}{\sqrt{2 m}}\sum_{k=1}^{
%m}((y_k/\e)^2-1), \psi_{\e,H}=\one_{\{t_{\e,m}>H\}}.
%
%Then, the following statement is true.
%
%and the hypothesis testing problem \nref{1.4} with the set under the
%alternative
%given by \nref{ell} with $q=2$. Let $a_k = k^\alpha$ and $\sigma_k = k^
%determined by \nref{MU}, and assume that $u_\e\to\infty$. Then, the
%family of tests $\psi_{\e, H}$, given by \nref{test2} with $m=\tilde
%m$ satisfying \nref{M} and $H=H_\e\to\infty$, is asymptotically
%minimax consistent, i.e., $ \a_{\e}(\psi_{\e,H_{\e}})\to0$ and
%there exists $c>0$ such that $\b_\e(\psi_{\e,H_\e},\Theta_\e)\to0$
%as $H_\e<(c+o(1))u_\e$.
%
%

Unlike the case $q=2$, the ``sparse'' case $q \in(0,2)$ is not
directly linked to Theorem~\ref{testing}; it will be considered
separately in Section~\ref{Smildlyq}.

%Also, the corresponding results for the Besov classes of functions
%will be presented in Section~\ref{Besq}.

%s4.3 #&#
\subsection{Application to severely ill-posed inverse problems with the
class of
analytic functions} \label{R2}

%th4.3 #&#
%
\begin{theorem}\label{T2}
%Consider the GSM \nref{ffgsm} and the hypothesis testing problem
Let $q\in(0,2], a_k = \exp(\alpha k)$ and $\sigma_k = \exp(\beta
k)$, $k\in\N$, $\alpha>0, \beta>0$.

\begin{longlist}[(a)]
\item[(a)]
The asymptotically minimax consistent family of tests
$\psi_{\e,H}$ are determined by the family of test statistics
$t_{\e}$ given by~(\ref{test1}) with $w_k$, $k\in\N$,
as in~(\ref{seq1}), and with
$\{\tilde{\eta}_k\}_{k\in\N}$ given by (11.17) with $m$
determined by (11.18) in~\cite{ISS12}.

\item[(b)] The separation rates are of the form
%
%e4.4 #&#
%
\begin{equation}
\label{seprateSevAnalF} r_\varepsilon^* = \varepsilon^{\alpha/(\alpha
+\beta)}.
\end{equation}
\end{longlist}
\end{theorem}

The proof is given in the supplementary material~\cite{ISS12}, Section 11.5.

%results obtained in Theorem~\ref{T2} hold true uniformly over
%$(\a,\b)\in\Sigma$ for any compact set $\Sigma\subset(0,

%re4.4 #&#
%
\begin{remark}
\label{w0FF1} We do not consider sharp asymptotics in this
case, since the assumption $w_0=o(1)$
does not hold for $\b>0$ in the case $q=2$. %Indeed,
%$$
%w_0=\frac{\max_{1\leq k \leq m}\tilde
%$$
\end{remark}

% statements for the related (well-posed or ill-posed) inverse problems
%are true:
%
%(a) if $\alpha\leq0$ and $\beta\geq0$, then $\gamma_{\e}(\Theta_{
%minimax testing is impossible.
%
%(b) If
%$\alpha>0$ and $\beta=0$, then the separation rates are of the form
%$
%r_\varepsilon^* =\e\log^{1/4}(\e^{-1}).
%$
%
%(c) If $\beta<0$, then we arrive at the classical separation rates,
%i.e.,
%$r_\e^*=\varepsilon$.
%
%The point (a) follows from Proposition 4.7 of Section 4.4.3
%in
%4.4.3 in~\cite{IS02}. The point (c) follows from Theorem 3.1 of
%Section 3.4.3 in~\cite{IS02}.}

%re4.5 #&#
%
\begin{remark}
\label{RemequivSevAn} Similar to Remark~\ref{RemequivSS}, the
asymptotics (11.19) in~\cite{ISS12} hold true for the sequences $a_k
\asymp\exp(\alpha k)$ and $\sigma_k \asymp\exp (\beta k)$, $k\in\N$,
$\alpha>0$, $\beta>0$. In this case, the separation rates are still of
the form~(\ref{seprateSevAnalF}).
%The Remarks~\ref{CompSetFF},
Remark~\ref{w0FF1} still applies to these cases, too.
\end{remark}

%Let us consider the case $q\in(0,2)$. Observe the embedding
%which yields $\g(\psi, \Theta_q(r_\e))\le\g(\psi, \Theta_2(r_\e))$.
%Therefore it suffices to establish the lower bounds.
%
%Take $m=\max\{k: r_\e\exp(\a k)\le1\}$, and consider the vector
%$\eta_m$ that contains only one non-zero coordinate
%$z_n=r_\e\exp(-\b n)$ at the place $m$. One easily check that
%$\eta_m\in\Theta_q(r_\e)$ for any $q>0$. Therefore one cannot
%distinguish between $H_0$ and $H_1$ if $z_n=o(\e)$, which is
%equivalent to $r_\e=o(r^*_\e)$ with $r^*_\e$ defined by
%true uniformly in $q\in(0,2)$ as well.

%s4.4 #&#
\subsection{Application to severely ill-posed inverse problems with
the Sobolev
class of functions}\label{R3}

%th4.4 #&#
%
\begin{theorem}\label{T3} %Consider the GSM \nref{ffgsm} and
%the hypothesis testing problem \nref{1.4}-\nref{ell}.
Let $q\in
(0,2], a_k = k^\alpha$ and $\sigma_k = \exp(\beta k)$, $k\in\N$,
$\alpha>0$, \mbox{$\beta>0$}.

\begin{longlist}[(a)]
\item[(a)]
The asymptotically minimax consistent family of tests
$\psi_{\e,H}$ are determined by the family of test statistics
$t_{\e}$ given by~(\ref{test1}) with $w_k$, $k\in\N$,
as in~(\ref{seq1}), and with
$\{\tilde{\eta}_k\}_{k\in\N}$ given by (11.20) with $m$
determined by (11.21) in~\cite{ISS12}.

\item[(b)] The separation rates are of the form
%
%e4.5 #&#
%
\begin{equation}
\label{seprateSevSobolF} r_\varepsilon^* = \bigl( \bigl(\log\bigl(
\varepsilon^{-1}\bigr) %-\alpha\log\log(\varepsilon^{-1})
\bigr)/\beta\bigr)^{-\alpha}.
\end{equation}
\end{longlist}
\end{theorem}

The proof is given in the supplementary material~\cite{ISS12}, Section 11.6.

%re4.6 #&#
%
\begin{remark}
\label{sharpsepratFFF}
A stronger result is possible in this case. In view
of (11.21) in~\cite{ISS12}, the relation
(\ref{seprateSevSobolF})
determines \textit{sharp separation rates}\vadjust{\goodbreak} $r_\varepsilon^*$ in the
following sense:

\begin{longlist}[(a)]
\item[(a)]
if $\liminf(r_\varepsilon/ r_\varepsilon^* ) >1$, then
$u_\varepsilon\to\infty$, that is,
$\gamma_\varepsilon(r_\varepsilon)\to0$;

\item[(b)] if $\limsup(r_\varepsilon/ r_\varepsilon^* )<1$, then
$u_\varepsilon\to0$, that is, $\gamma_\varepsilon(r_\varepsilon
)\to
1$, and the minimax testing is impossible.
\end{longlist}

Moreover the relation $
(r_\varepsilon^*)^{-1/\a}= ( (\log(\varepsilon^{-1})-\alpha
\log\log
(\varepsilon^{-1}) )
/\beta)+O(1),
$
determines the sharper separation rates $r_\varepsilon^*$ in the
following sense:

\begin{longlist}[(c)]
\item[(c)]
if
$\lim\inf(r_\varepsilon^{-1/\a}-({r_\varepsilon^*})^{-1/\a
} )=-\infty$,
then $u_\varepsilon\to\infty$, that is,
$\gamma_\varepsilon(r_\varepsilon)\to0$;

\item[(d)] if
$\lim\sup(r_\varepsilon^{-1/\a}-({r_\varepsilon^*})^{-1/\a
} )=+\infty$,
then $u_\varepsilon\to0$, that is,
$\gamma_\varepsilon(r_\varepsilon)\to1$, and the testing is
asymptotically impossible.
\end{longlist}
\end{remark}

%{\rm It is also easy to see that the asymptotic results in Theorem
%uniformly over
%$(\a,\b)\in\Sigma$ for any compact set $\Sigma\subset(0, \infty)

%re4.7 #&#
%
\begin{remark}
\label{w0FF2} We do not consider sharp asymptotics in this
case, since the assumption $w_0=o(1)$ does not hold for $\b>0$ in
the case $q=2$. %Indeed,
%$$
%B_1>0; \sqrt{2\sum_{1\le k\le m}\tilde
%$$
%and, therefore,
%$$
%w_0=\frac{\max_{1\leq k \leq m}\tilde\eta_k^2}{\sqrt{2\sum_{1\le
%k\le m}\tilde\eta_k^4}}\ge B\frac{z_0^2 m \exp(2\beta m)}{z_0^2
%m\exp(2\beta m)}\asymp1 \not\to0, B>0.
%$$
\end{remark}

%(well-posed or ill-posed) inverse problems are true:
%
%(a) if $\alpha\leq0$ and $\beta\geq0$, then this case corresponds to
%severely ill-posed
%(well-posed if $\beta=0$) inverse problems with the class of analytic
%functions
%(see Remark~\ref{Rem2}) and, hence, $\gamma_{\e}(\Theta_{\e}) \equiv
%1$, i.e., minimax
%testing is impossible.
%
%(b) If $\alpha>0$ and $\beta=0$, then this case corresponds to a
%well-posed inverse problem with the Sobolev class of functions (see
%Remark~\ref{FFF1}) and, hence, the separation rates are of the form
%$
%r_\varepsilon^* = \varepsilon^{4 \alpha/(4\alpha+1)}.
%$
%(c) If $\beta<0$, then we arrive at the classical separation rates,
%i.e., $r_\e^*=\varepsilon$.
%
%in
%in Section~\ref{R1} and follows from Remark~\ref{FFF1}. The point
%(c) follows from Theorem 3.1 of Section 3.4.3 in~\cite{IS02}.}
%

%re4.8 #&#
%
\begin{remark}
\label{equivR3FF2} Similar to Remark~\ref{RemequivSS}, the
asymptotics (11.22) in~\cite{ISS12} and the sharp separation rates
(\ref{seprateSevSobolF}) mentioned in Remark~\ref{sharpsepratFFF}
hold true for the sequences $\sigma_k \asymp\exp(\beta k)$,
$k\in\N$, $\beta>0$. The dependence on the sequence
$\{a_k\}_{k\in\N}$ is, however, more delicate. One can actually show
that the sharp separation rates~(\ref{seprateSevSobolF}) mentioned
in Remark~\ref{sharpsepratFFF} are still of the same form for $a_k
\sim k^{\alpha}$, $k\in\N$, $\alpha>0$.
%The Remarks~\ref{CSFF2},
Remark~\ref{w0FF2} still applies to these cases too.
\end{remark}

%s4.5 #&#
\subsection{Application to mildly ill-posed inverse problems with the
class of
analytic functions}\label{R4}

Here, we consider the ``standard'' case $q =2$. [The ``sparse'' case $q
\in(0, 2)$
will be discussed in Remark~\ref{remsparsecase}.]

%th4.5 #&#
%
\begin{theorem}\label{T4}
%Consider the GSM \nref{ffgsm} and the hypothesis testing problem
Let $q=2, a_k = \exp(\alpha k)$
and $\sigma_k = k^\beta$, $k\in\N$, $\alpha>0, \beta>0$.

\begin{longlist}[(a)]
\item[(a)]
The sharp asymptotics~(\ref{G}) hold with the the value $u_\e$ of
the extreme problem~(\ref{D1}) determined by (11.29) in~\cite{ISS12}.

\item[(b)] The asymptotically
minimax family of tests are determined by the test statistics $t_{\e}$
given by~(\ref{test1})
with $w_k$, $k\in\N$ and $w_0$ as in~(\ref{seq1}), and with
$\{\tilde{\eta}_k\}_{k\in\N}$ given
by (11.24) with $m$ determined by (11.30) in~\cite{ISS12}.

\item[(c)] The separation rates are of the form
%
%e4.6 #&#
%
\begin{equation}
\label{seprateAnSovF} r_\e^*=\e\bigl(\log\e^{-1}
\bigr)^{\beta+1/4}.
\end{equation}
\end{longlist}
\end{theorem}

The proof is given in the supplementary material~\cite{ISS12}, Section 11.7.

%{\rm Under the assumption \nref{sigma}, it is immediately seen that
%the the asymptotic
%results in Theorem~\ref{T4} hold true for any $\beta> -1/4$. Hence,
%the separation rates
%consideration but also for a range of corresponding well-posed problems
%($\beta\in(-1/4, 0]$).}

%re4.9 #&#
%
\begin{remark}
\label{nFF1F}
%It is also easy to see that,
%the asymptotic results in Theorem~\ref{T4} hold true
%uniformly over
%$(\a,\b)\in\Sigma$ for any compact set $\Sigma\subset(0, \infty)
%(-1/4,\infty)$. In particular,
%uniformly over $(\a,\b)\in\Sigma$, the efficient
%dimension $m=m_\e(\a,\b)$ satisfies
%m_\e(\a,\b)\sim\frac{2\log(\e^{-1})-\log(u_\e)}{2\a}\asymp
%
%%%%% Rewritten Version for the revision %%%%%
It is also easy to see that, uniformly over $(\a,\b)\in\Sigma$,
for any compact set $\Sigma\subset\R_+^2$,
the efficient dimension $m=m_\e(\a,\b)$ satisfies
%
%e4.7 #&#
%
\begin{eqnarray}
\label{m**} m_\e(\a,\b)&\sim&\frac{2\log(\e^{-1})-\log(u_\e)}{2\a}\nonumber\\[-8pt]\\[-8pt]
&\asymp&
\log\bigl(
\e^{-1}\bigr) \qquad\mbox{as } \log(u_\e)=o\bigl(\log\bigl(
\e^{-1}\bigr)\bigr).\nonumber
\end{eqnarray}
\end{remark}

%{\rm Using Theorem 3.1 of Section 3.4.3 in~\cite{IS02}, and the
%results derived in
%Sections~\ref{R1} and~\ref{R2},
%we can easily see that the following statements for the related
%(well-posed or ill-posed)
%inverse problems are true:
%
%(a) if $\alpha\leq0$ and $\beta\geq-1/4$, then
%$\gamma_{\e}(\Theta_{\e}) = 1$, i.e., minimax testing is impossible.
%
%(b) If $\alpha>0$ and $\beta=-1/4$, then
%$$
%w_0\le\frac{Bz_0^2}{z_0^2\log m}\to0,
%m\sim\log(r_\e^{-1/\alpha}),
%$$
%and on scrutinizing the arguments to arrive at Theorem~\ref{T4}, one
%can easily see that we get the sharp asymptotics \nref{G} with $
%u_\e^2\sim r_\e^4 \e^{-4}(2\log\log(r_\e^{-1/\alpha}))^{-1}, $ and,
%in view of \nref{m**}, the separation rates are of the form $
%r_\e^*=\e(\log\log\e^{-1})^{1/4}.
%$
%
%(c) If $\beta<-1/4$, then we arrive at the classical separation rates,
%i.e.,
%$r_\e^*=\varepsilon$.}

%re4.10 #&#
%
\begin{remark}\label{R10}
As in Remark~\ref{RemequivSS}, the asymptotics (11.30) in
\cite{ISS12} hold true for the sequences $a_k \sim\exp(\alpha k)$ and
$\sigma_k \sim k^{\beta}$, $k\in\N$, $\alpha>0$, $\beta>0$. Similar
rate asymptotics hold true for the sequences $a_k \asymp\exp(\alpha
k)$ and $\sigma_k \asymp k^{\beta}$, $k\in\N$, $\alpha>0$, $\beta>0$.
In both cases, the separation rates are still of the form~(\ref{seprateAnSovF}).
%The Remarks~\ref{FFF2},~\ref{nFF1F} and~\ref{nFF2F} still apply to
%these cases too.
\end{remark}

%re4.11 #&#
%
\begin{remark}
\label{remsparsecase}
The rate asymptotics of Theorem~\ref{T4} hold true uniformly in $q
\in[\delta,2]$, for any
$\delta\in(0,2)$. Indeed, in view of the embedding
%
%e4.8 #&#
%
\begin{equation}
\label{embedd} \Theta_q\subset\Theta_2,\qquad
\Theta_q(r_\e)\subset\Theta_2(r_\e)
\end{equation}
for the sets defined by~(\ref{ell}), it suffices to establish the
lower bounds. Let $ \Theta_q^\a(L)=\{\eta\in l^2\dvtx
\sum_{k\in\N}|\exp(\a k)k^\b\eta_k|^q\le L^q\}. $ [Note that the
set $\Theta_q^\a(1)$ just corresponds to the set under the
alternative considered above.] For any $\a_1=\a+\delta, \delta>0$,
we have the embedding
%
%e4.9 #&#
%
\begin{equation}
\label{embedd2}\qquad \Theta_2^{\a_1}(cL)\subset
\Theta_q^\a(L),\qquad c=c(q,\delta)=\bigl(\exp\bigl(2q
\delta/(2-q)\bigr)-1\bigr)^{(2-q)/2q}
\end{equation}
and $c(q,\delta)\to\exp(\delta)$ as $q\nearrow2$. Using H\"older's
inequality, the above embedding follows easily on noting that
\begin{eqnarray*}
&&\sum_{k\in\N}\bigl|e^{\a k}k^\b
\eta_k\bigr|^q \le\biggl(\sum_{k\in\N
}
\bigl(e^{\a
k}k^\b\eta_k\bigr)^2
\biggr)^{{2}/{q}} \biggl(\sum_{k\in\N}
e^{-2kq\delta/(2-q)} \biggr)^{1-{q}/{2}}
\\
&&\qquad= \biggl(c^{-2}\sum_{k\in\N}
\bigl(e^{\a k}k^\b\eta_k\bigr)^2
\biggr)^{{2}/{q}}.
\end{eqnarray*}
%
%k)k^\b\eta_k|^q\exp(-q\delta)\\
%&\le& \left(\sum_{k\in\N}(\exp(\a k)k^\b\eta_k)^2\r)^{2/q}
%&=&\left(c^{-2}\sum_{k\in\N}(\exp(\a k)k^\b\eta_k)^2\r)^{2/q}.
Since the separation rates from Theorem~\ref{T4} do not depend on
$\a$ and $c$, in view of Remark~\ref{R10}, the rate asymptotics of
Theorem~\ref{T4} hold true uniformly in $q \in[\delta,2]$, for any
$\delta\in(0,2)$.
\end{remark}

\subsection{Application to extremely ill-posed inverse problems with
the class of
generalized analytic functions}\label{R5}

We consider the case $q=2$ only. Assume that $\{a_k\}_{k\in\N}$ and
$\{\sigma_k\}_{k\in\N}$ are increasing sequences such that
%
%e4.10 #&#
%
\begin{equation}\quad
\label{assump} \lim_{k\to\infty}\sigma_{k+1}/\sigma_k\to
\infty,\qquad \liminf_{k\to\infty}a_{k+1}/a_k=c,\qquad c\in(1,\infty].
\end{equation}
In order to describe the asymptotics of the value $u_\e=u_\e(r)$ of the
extreme problem
(\ref{D1}), we introduce the following functions.

Let $m\in\N$, $m \geq2$, $\Delta^*_m=[1/a_{m},
1/a_{m-1}]$, and for $r<1/a_1$ take $m=m(r)\geq2$ such that $r \in
\Delta^*_m$. Consider
now the piecewise quadratic (in $r^2$) function defined by
%
%e4.11 #&#
%
\begin{equation}\qquad
\label{e1FFF} \bigl(u_\e^*(r)\bigr)^2=\frac{1}{2\e^4(a_m^2-a_{m-1}^2)^2}
\biggl(\frac{(a_m^2
r^2-1)^2}{\sigma_{m-1}^4}+\frac{(1-a_{m-1}^2
r^2)^2}{\sigma_{m}^4} \biggr),
\end{equation}
and the piecewise linear (in $r^2$) function defined by
%
%e4.12 #&#
%
\begin{equation}
\label{elinFFF} u_\e^{\mathrm{lin}}(r)=\frac{1}{\e^2(a_m^2-a_{m-1}^2)} \biggl(
\frac{a_m^2
r^2-1}{\sigma_{m-1}^2}+\frac{1-a_{m-1}^2
r^2}{\sigma_{m}^2} \biggr).\vadjust{\goodbreak}
\end{equation}

%Then, the following statement is true.

%th4.6 #&#
%
\begin{theorem}\label{TheorFFF} Let $u_\e=u_\e(r)$ be the value of the
extreme problem
(\ref{D1}). Let $(u_\e^*(r))^2$ be the
piecewise quadratic (in $r^2$) function defined by~(\ref{e1FFF}), and let
$u_\e^{\mathrm{lin}}(r)$ be the
piecewise linear (in $r^2$) function defined by~(\ref{elinFFF}), where
$\{a_k\}_{k\in\N}$
and $\{\sigma_k\}_{k\in\N}$ be increasing sequences satisfying (\ref
{assump}).

\begin{longlist}[(a)]
\item[(a)]
(Sharp asymptotics of $u_\e$.) The families $u_{\e}(r)$ and $u_{\e
}^*(r)$ are related~by
%
%e4.13 #&#
%
\begin{equation}
\label{e2} u_\e(r_\e)\sim u_\e^*(r_\e)
\qquad\mbox{as } r_{\e}\to0.
\end{equation}

\item[(b)] (Rate asymptotics of $u_\e$.) The families $u_{\e}(r)$ and $u_{\e
}^{\mathrm{lin}}(r)$ are related~by
%
%e4.14 #&#
%
\begin{equation}
\label{elin1} \hspace*{32pt} u_\e^{\mathrm{lin}}(r_\e) \bigl(1/2+o(1)
\bigr)\le u_\e(r_\e)\le u_\e^{\mathrm{lin}}(r_\e)
\bigl(1/\sqrt{2}+o(1)\bigr) \qquad\mbox{as } r_\e\to0.
\end{equation}

\item[(c)] (Distinguishability conditions.) Consider the GSM~(\ref{ffgsm}) and
the hypothesis testing problem~(\ref{14}) and~(\ref{ell}). Then
\begin{eqnarray*}
\g_\e(r_\e)&\to&0 \quad\mbox{if and only if}\quad
u_\e^{\mathrm{lin}}(r_\e)\to\infty;\\
\g_\e(r_\e)&\to&1 \quad\mbox{if and only if}\quad
u_\e^{\mathrm{lin}}(r_\e)\to0.
\end{eqnarray*}
%
%(d) (separation rates) The separation rates $r_{\e}^*$ are
%determined by the relation $u_{\e}^{lin}(\Theta(r_{\e}^*))\asymp1$.
\end{longlist}
\end{theorem}

The proof is given in the supplementary material~\cite{ISS12}, Section 11.9.

%re4.12 #&#
%
\begin{remark}
It is easy to see that relation~(\ref{elin1}) is true uniformly over all
sequences
$\{a_k\}_{k\in\N}$ and $\{\sigma_k\}_{k\in\N}$
such that $\sigma_{k+1}/\sigma_k \geq B_k$, $B_k \to\infty$, and
$a_{k+1}/a_k > c$,
as $k\geq k_0$, $k_0 \geq1$.
\end{remark}

%re4.13 #&#
%
\begin{remark}
\label{remlast2val}
We do not consider sharp asymptotics in
this case, since the assumption $w_0=o(1)$ does not hold under
assumption~(\ref{assump}). %Indeed, using
%that
%$
%$
%and, therefore,
%$$
%w_0=\frac{\max_{1\leq k \leq m}\tilde
%$$
\end{remark}

%re4.14 #&#
%
\begin{remark}
\label{remsharp}
The relation $u_\e^{\mathrm{lin}}(r_\e^*)\asymp1$
determines the separation rates~$r_\e^*$ that are rather sharp in the follows
sense. Let $r_\e^*=a_m^{-1}$ for some $m\in\N,\break m\to\infty$, and let
$r^2=(1+b)(r_\e^*)^2\in(a_m^{-2},a_{m-1}^{-2}), b>0$. Then, one has
$
u_\e^{\mathrm{lin}}(r)=u_\e^{\mathrm{lin}}(r_\e^*)(1+k_mb),
$
where, as $m\to\infty$,
\begin{eqnarray*}
k_m&=&\frac{\sigma_m^2}{\sigma_{m-1}^2}\frac{1-(\sigma
_{m-1}a_{m-1}/\sigma_ma_m)^2} {
1-(a_{m-1}/a_m)^2}
\\
&\sim&\frac{\sigma_m^2}{\sigma_{m-1}^2(1-(a_{m-1}/a_m)^2)}\asymp\frac
{\sigma_m^2}{\sigma_{m-1}^2}\to\infty.
\end{eqnarray*}
Therefore, in order to obtain $u_\e^{\mathrm{lin}}(r_\e)\to\infty$, it suffices
to take $r_\e=r_\e^*(1+\delta)$ for any $\delta>0$. On the other
hand, let $r^2=(1-b)(r_\e^*)^2\in(a_{m+1}^{-2},a_m^{-2}), b\in
(0,1)$. Then, similarly, one has
$
u_\e^{\mathrm{lin}}(r)=u_\e^{\mathrm{lin}}(r_\e^*)(1-l_mb),
$
where, as $m\to\infty$,
\[
l_m=\frac{1-(\sigma_ma_m/\sigma_{m+1}a_{m+1})^2} {1-(a_m/a_{m+1})^2}
\sim\frac{1}{1-(a_m/a_{m+1})^2}\asymp1.
\]
If $a_{m+1}/a_m\to\infty$ as $m\to\infty$, then, in order to obtain
$u_\e^{\mathrm{lin}}(r_\e)\to0$, one~needs to take $r_\e$ such that
$r_\e/r_\e^*\to0$, and if $a_{m+1}/a_m\to c, 1<c<\infty$, then,~for
\mbox{$u_\e^{\mathrm{lin}}(r_\e)\to0$}, one needs to take $r_\e<r_\e^*/c$.
Thus, the conditions for distinguish\-ability and
nondistinguishability could be nonsymmetric in these
problems.\vadjust{\goodbreak}\looseness=-1
\end{remark}

%re4.15 #&#
%
\begin{remark}
\label{remexample} Let us consider the example
\[
a_k=\exp\bigl(\a k^\tau\bigr),\qquad \alpha> 0, \tau\ge1,\qquad
\sigma_k=\exp\bigl(\b k^\gamma\bigr),\qquad \beta>0, \gamma>1.
\]
For the moment, let us forget that $m\in\N$ and define $m=m(r)\in\R_+$
by the
equality $ r=a_m^{-1}$, that is, $m(r)=(\log(r^{-1})/\a)^{1/\tau}$. Set
also
\[
\hat u_\e(r)=(\e a_{m(r)}\sigma_{m(r)})^{-2}=(r/
\e)^2\exp\bigl(-2\b\bigl(\log\bigl(r^{-1}\bigr)/\a
\bigr)^{\gamma
/\tau}\bigr).
\]
Observe that $\hat u_\e(r)= u_\e^{\mathrm{lin}}(r)$ as $r=a_m^{-1}, m\in
\N$. On the other hand, one can check that the function $\hat
u_\e(r)$ is a convex function in $r^2$ for $r>0$ small enough.
Therefore, $\hat u_\e(r)< u_\e^{\mathrm{lin}}(r)$ as $r\not=a_m^{-1}$ for any
$m\in\N$, and the condition $\hat u_\e(r_\e)\to\infty$ implies
$u_\e^{\mathrm{lin}}(r_\e)\to\infty$. However, it is possible that
$u_\e^{\mathrm{lin}}(r_\e)\to\infty$ when $\hat u_\e(r_\e)=O(1)$, in general.
For instance, let $\tau=\gamma$. Then
\[
\hat u_\e(r)=\e^{-2}r^{2+2\b/\a}.
\]
If $r_\e=a_m^{-1}$ and
$u_\e^{\mathrm{lin}}(r_\e)\asymp1$, then it was noted in Remark
\ref{remsharp} that $u_\e^{\mathrm{lin}}(r_\e(1+\delta))\to\infty$ for any
$\delta>0$, but $\hat u_\e(r_\e(1+\delta))\asymp1$. The same holds
for $\gamma<\tau$.

The relation $\hat u_\e(\tilde r_\e)\asymp1$ determines the family
$ \tilde r_\e$. %(compare with the `minimax rates of testing' in
Note that if $r_\e/\tilde r_\e\to\infty$, then $\hat
u_\e(r_\e)\to\infty$, and since $u_\e^{\mathrm{lin}}(r)\ge\hat u_\e(r)$,
this yields $u_\e^{\mathrm{lin}}(r_\e)\to\infty$ and $\g_\e(r_\e)\to0$ by
Theorem~\ref{TheorFFF}. However, this family is not a family of
separation rates, at least if $\gamma\le\tau$, because the
condition $r_\e/\tilde r_\e\to0$ does not guaranty that
$\g_\e(r_\e)\to1$.

More precisely, there exists a sequence $\e_m\to0$ and $\hat
r_{m}=o(\tilde r_{\e_m})$ such that $\g_{\e_m}(\hat r_{m})\to0$.
In fact, observe that if $\gamma\le\tau$, then the function $\hat
u_\e(r)$ satisfies [uniformly over $\e>0$ since $\e^2$ is a factor
in $u_\e(r)$]
\[
\hat u_\e(Br)\asymp\hat u_\e(r) \quad\mbox{iff}\quad B\asymp1,\qquad
r\to0.
\]
Take a sequence $m\to\infty$ and put $r_m=a_m^{-1},
(r_m^{(1)})^2=r_m^2(1+\delta_m), \hat r_m^2=r_m^2(1+\delta)$, where
$\delta_m\to0, \delta_m\sigma_m^2/\sigma_{m-1}^2\to\infty,
\delta>0$. Observe that similarly to evaluations in Remark
\ref{remsharp}, one has, uniformly over $\e>0$,
\[
\hat u_\e(\hat r_m)\asymp\hat u_\e\bigl(
r_m^{(1)}\bigr)\asymp\hat u_\e(
r_m)= u_\e^{\mathrm{lin}}( r_m)\ll
u_\e^{\mathrm{lin}}\bigl( r_m^{(1)}\bigr)\ll
u_\e^{\mathrm{lin}}( \hat r_m),\qquad m\to\infty.
\]
Take now $\tilde r_m$ and $\e_m$ such that
$
\hat u_{\e_m}(\tilde r_m)\asymp u_{\e_m}^{\mathrm{lin}}( r_m^{(1)})\asymp1.
$
This implies\break $\hat u_{\e_m}(\tilde r_m)\gg\hat u_{\e_m}( r_m)$ and
$\tilde r_m\gg r_m\asymp\hat r_m$. By construction, we see that the
sequence $\tilde r_m$ satisfies $\hat u_{\e_m}(\tilde r_m)\asymp1$
and $\hat r_m=o(\tilde r_m)$, but $u_{\e_m}^{\mathrm{lin}}(\hat
r_m)\to\infty$, which yields $\g_{\e_m}(\hat r_{m})\to0$.
%
% This remark show that $\tilde r_\e$ are not the separation rates
%in the sense of \nref{rates1}.
\end{remark}

\subsection{\texorpdfstring{Mildly ill-posed inverse problems with $l^q$-ellipsoids
for Sobolev
classes of functions: The ``sparse'' case $q \in(0,2)$}
{Mildly ill-posed inverse problems with $l^q$-ellipsoids
for Sobolev
classes of functions: The ``sparse'' case q in (0,2)}}\label{Smildlyq}

%Here, we consider mildly ill-posed inverse problems with
%$l^q$-ellipsoids for the Sobolev classes of functions, for the
%``sparse'' case $q \in(0,2)$.
Unlike the ``standard'' case $q=2$, the sharp and rate
optimality results for the ``sparse'' case $q \in(0,2)$ are of
different nature and are not directly linked with Theorem
\ref{testing}, but can be obtained from a hitherto unknown link with
results obtained in another context and presented in Sections
4.4.2--4.4.3
%and Sections 6.2.1-6.2.2
of~\cite{IS02}.
%Since these results are scattered in the cited reference, and not not
%immediately seen as ill-posed inverse problems.
For completeness and an immediate access to these results, we formulate
and present them below.\vadjust{\goodbreak}

Consider
the extreme problem
%
%e4.15 #&#
%
\begin{equation}
\label{extrel} u_\e^2= 2\inf\sum
_{i\in\N} h_i^2\sinh^2
\bigl(z_i^2/2\bigr),
\end{equation}
where the infimum is taken over sequences $(h_i, z_i), h_i\in[0,1],
z_i\ge0, i\in\N$, such that
%
%e4.16 #&#
%
\begin{equation}
\label{extrelconstr} \sum_{i\in\N}
i^{2\b}h_iz_i^2\geq\bigl(\tilde
r_\e^2/\e^2\bigr),\qquad \sum
_{i\in\N} i^{q(\a+\b
)}h_iz_i^q
\leq\bigl(1/\e^q\bigr),
\end{equation}
where $\tilde r_\e=r_\e(1-\delta_\e), \delta_\e>0, \delta_\e
\to
0, \delta_\e\log(\e^{-1})\to\infty$.

Set $\lambda=(\alpha+\beta)/2-\beta/q$. If $\lambda>0$, then there
exist extreme sequences $h_{i,\e}\in
(0,1], z_{i,\e}>0$, in the problem~(\ref{extrel}) and (\ref
{extrelconstr}), and we have the asymptotics of the form
%
%e4.17 #&#
%
\begin{equation}
\label{u1ell} u_\e^2\sim
c_0nh_0^2,
\end{equation}
where the quantities $n=n_\e$ and $h_0=h_{0,\e}$
are determined by the relations
%
%e4.18 #&#
%
\begin{equation}
\label{u2el} c_1 n^{\b+1/2}h_0^{1/2}
\sim r_\e/\e c_2 n^{\a+\b+1/q}h_0^{1/q}
\sim1/\e
\end{equation}
for some constants $c_l=c_l(\a,\b,q)>0, l=0,1,2$, which, in turn, imply
%
%e4.19 #&#
%
\begin{equation}\quad
\label{extru} u_\e\sim c_3\e^{-(2\a+1/q-1/2)/(\a+\b(1-2/q))}r_\e^{(2(\a
+\b)+1/q)/(\a
+\b(1-2/q))}
\end{equation}
for some constant $c_3=c_3(\a,\b,q)>0$. (The quantity $ n=n_\e\to
\infty
$ plays the role of the ``efficient dimension'' in the problem.)

Set
$
Q_{\e,i}=\sqrt{2(\log i+\log\log i+2\log\log(\e^{-1}))}
$
and consider the events
%
%e4.20 #&#
%
\begin{equation}
\label{thresh} {\CY}_\e=\Bigl\{y=\{y_i
\}_{i\in\N}\dvtx \Bigl(\sup_{i\in\N}|y_i|\big/(\e
Q_{\e,i})\Bigr)>1\Bigr\}
\end{equation}
and the following families of test statistics:
%
%e4.21 #&#
%
\begin{equation}\quad
\label{0stat} l_\e(y)=u_\e^{-1}\sum
_{i\in\N} h_{\e,i}\xi(y_i/
\e,z_{\e
,i}),\qquad \xi(t,z)=e^{z^2/2}\cosh(tz)-1
\end{equation}
and tests
%
%e4.22 #&#
%
\begin{equation}
\label{testGD} \psi^{G}_{\e,H}=
\one_{\{l_\e(y)>H\}\cap{\CY}_\e},\qquad
\psi^{D}_{\e}=\one_{{\CY}_\e},\qquad
\psi^{D}_{\e,\a}=\a+(1-\a)\one_{{\CY}_\e}.\hspace*{-35pt}
\end{equation}

%Then, the following statement is true.

%th4.7 #&#
%
\begin{theorem}\label{T1Sparse}
%Consider the GSM \nref{ffgsm} and the hypothesis testing problem
Let $q\in(0,2), a_k = k^\alpha$ and
$\sigma_k = k^\beta$, $k\in\N$, $\alpha>0, \beta>0$, and set
$\lambda=(\alpha+\beta)/2-\beta/q$.

\begin{longlist}[(a)]
\item[(a)]
If $\lambda>0$, then the sharp asymptotics are of the Gaussian
type~(\ref{G}) with $u_\e$ from~(\ref{extrel}). The tests
$\psi^{G}_{\e,H}$ of the form~(\ref{testGD}) with $H=H^{(\a)}$ and
$H=u_\e/2$ are asymptotically minimax, that is,
\begin{eqnarray*}
\a_{\e}\bigl(\psi_{\e,H^{(\a)}}^{G}\bigr) &\leq&\a+o(1),\qquad
\b_{\e}\bigl(\Theta(r_{\e}),\psi_{\e,H^{(\a)}}^{G}
\bigr) = \b_{\e}(r_{\e},\a)+o(1),
\\
\g_{\e}\bigl(\Theta(r_{\e}), \psi_{\e,u_\e/2}^{G}
\bigr) &=& \g_{\e}(r_{\e})+o(1).
\end{eqnarray*}

\item[(b)] If $\lambda\le0$, then the sharp asymptotics are of the
following degenerate type: $ \b_\e(r_\e,\a) =
(1-\a)\Phi(-D_\e)+o(1), \g_\e(r_\e) =\Phi(-D_\e)+o(1), $
where
$
D_\e=n_\e^{-\b} r_{\e}/{\e}-\sqrt{2\log(n_\e)}, n_\e=r_\e^{-1/\a}.
$
The tests~$\psi^{D}_{\e}$ and (the randomized) tests~$\psi^{D}_{\e,\a}$
of the form~(\ref{testGD}) are asymptotically minimax, that is,
\begin{eqnarray*}
\a_\e\bigl(\psi^{D}_{\e,\a}\bigr) &=& \a+o(1),\qquad
\b_{\e}\bigl(\Theta(r_{\e}),\psi_{\e,\a}^{D}
\bigr) = \b_{\e}(r_{\e},\a)+o(1),
\\
\g_{\e}\bigl(\Theta(r_{\e}), \psi_{\e}^{D}
\bigr) &=&\g_{\e}(r_{\e})+o(1).
\end{eqnarray*}

\item[(c)] If $\lambda>0$, then the separation rates are of the form
\[
r_\varepsilon^*= \varepsilon^{(2\alpha+1/q-1/2))/(2(\alpha+\beta)+1/q)}.
\]

\item[(d)] If $\lambda\le0$, then the sharp separation rates are of the
form
\[
r^*_\e=\Lambda\e^{\a/(\a+\b)}\bigl(\log\bigl(\e^{-1}
\bigr)\bigr)^{\a/2(\a+\b
)},\qquad \Lambda=\bigl(2/(\alpha+\b)\bigr)^{\a/2(\a+\b)}.
\]
\end{longlist}
\end{theorem}

Theorem~\ref{T1Sparse} is obtained by taking into account the
minimax hypothesis testing framework considered in Section
\ref{minimax} and Theorems 4.5 and 6.1 in~\cite{IS02}, noting
(from their proofs) that the events (thresholding rule) (4.155) in
\cite{IS02} can be replaced by the events (thresholding rule)
(\ref{thresh}). Its proof is omitted. The key ideas of the
study are discussed in the supplementary material~\cite{ISS12}, Section 10.
%(Note that there is a misprint in (4.155) in~\cite{IS02}:
%$\cup$ should be replaced by $\cap$.)

%The scheme of the proofs are given in the Appendix,
%Section~\ref{A2Sparse}.

%Using the general minimax hypothesis testing framework considered in
%Section~\ref{minimax}, Theorem~\ref{T1Sparse}
%% and relations \nref{u1_ell}, \nref{u2_el}
%can be seen as a particular cases of Theorem 4.5 in~\cite{IS02}, p.
%170, and Theorem
%6.1 in~\cite{IS02}, pp. 225-226.

%re4.16 #&#
%
\begin{remark}
\label{RemequivqSob}
Similar to Remark~\ref{RemequivSS}, we
get the sharp asymptotics~(\ref{extru}) for the sequences $a_k \sim
k^{\alpha}$ and $\sigma_k \sim k^{\beta}$, $k\in\N$, $\alpha>0,
\beta>0$, and similar rate asymptotics for the sequences $a_k \asymp
k^{\alpha}$, $\sigma\asymp k^{\beta}$, $k\in\N$, $\alpha>0,
\beta>0$. In both cases, the separation rates are still of the form
given in Theorem~\ref{T1Sparse}.
%The Remarks~\ref{FFF1},~\ref{FFF1FF} and~\ref{Rem1} still apply to
%these cases too.
%Remark~\ref{FFF1FF} still applies to these cases too.
\end{remark}

%s5 #&#
\section{Minimax signal detection in ill-posed inverse problems:
Adaptivity and rate optimality}\label{Ad}

The families of tests described in Section
\ref{subsecgenresF} (except those described in the supplementary
material~\cite{ISS12}, Theorem 7.1)
depend on a parameter $\kappa\in\Sigma\subset\R_+^{n} \times
(0,2]$, $n\geq
2$, associated with the sequences $\{a_k\}_{k\in\N}$ and
$\{\sigma_k\}_{k\in\N}$, and $q\in(0,2]$, that are involved in the
ill-posed inverse problems under consideration, that are usually
unknown in practice. For example, if $a_k = \exp(\a k^{\tau})$ and
$\sigma_k = \exp(\b k^{\gamma})$, $k\in\N$, $\a>0$, $\tau\geq1$,
$\b>0$,\break $\gamma> 1$, and if $q\in(0,2]$, then $\kappa\in
\Sigma=\{(\a,\tau,\b,\gamma, q)\} = (0, \infty) \cup[1, \infty)
\cup
(0,\infty) \cup(1, \infty) \cup(0,2]\subset\R_+^4 \times(0,2]$.

It is of paramount importance to construct families of tests that do not
depend on
the unknown parameter $\kappa$ and, at the same time, provide the best possible
asymptotical minimax efficiency.
These families of tests are called \textit{adaptive} (to the parameter
$\kappa$),
and the formal setting is as follows.

%s5.1 #&#
\subsection{Adaptive distinguisability and adaptive separation
rates}\label{adistasraF}

Let a set $\Sigma=\{\kappa\}$ and a family $r_\e(\kappa),
\kappa\in\Sigma$, be given, where $\e>0$ is small. Let the set
$\Theta_{\e}(\kappa, r_{\e}(\kappa))$ be determined by the
constraints~(\ref{ell}) with $a_k=a_k(\kappa),
\sigma_k=\sigma_k(\kappa)$, $k \in\N$, $q=q(\kappa)$,\vadjust{\goodbreak} and
$r_\e=r_{\e}(\kappa)$, and set
$
\Theta_{\e}(\Sigma)=\bigcup_{\kappa\in\Sigma}\Theta_\e(\kappa
,r_{\e
}(\kappa)).
$
We are interested in the following hypothesis testing problem:
\[
H_0\dvtx \eta=0\quad\mbox{versus}\quad H_1\dvtx \eta\in
\Theta_{\e}(\Sigma). %\end{cases}
\]
We are aiming to find conditions for either
$\g_\e(\Theta_{\e}(\Sigma))\to1$ or $\g_\e(\Theta_{\e}(\Sigma
))\to
0$, and to constructing asymptotically minimax adaptive consistent
families of tests
$\psi^{\mathrm{ad}}_\e$ such that
$\g_\e(\Theta_{\e}(\Sigma),\psi^{\mathrm{ad}}_\e)\to0$ as
$\g_\e(\Theta_{\e}(\Sigma))\to0$.

Let $u_\e(\kappa)=u_\e(\kappa, r_{\e}(\kappa))$ be the value of the
extreme problem
(\ref{D1}) for
the set $\Theta_\e=\Theta_{\e}(\kappa,r_{\e}(\kappa))$. Set
$
u_\e(\Sigma)=\inf_{\kappa\in\Sigma}u_\e(\kappa).
$
We are interested in how large $u_\e(\Sigma)$ should be in order to
provide the relation $\g_\e(\Theta_{\e}(\Sigma))\to0$.
We say that the family $u_\e^{\mathrm{ad}}=u_\e^{\mathrm{ad}}(\Sigma)\to\infty$
characterizes \textit{adaptive distinguishability} if there exist
constants $0<d=d(\Sigma)\le D=D(\Sigma)<\infty$ such that
\begin{eqnarray*}
\g_\e\bigl(\Theta_{\e}(\Sigma)\bigr)\to1 \qquad\mbox{as } \lim
\sup_{\kappa\in\Sigma}u_\e(\kappa)/u_\e^{\mathrm{ad}}&<&d,
\\
\g_\e\bigl(\Theta_{\e}(\Sigma)\bigr)\to0 \qquad\mbox{as } \lim
\inf_{\kappa\in\Sigma}u_\e(\kappa)/u_\e^{\mathrm{ad}}&>&D.
\end{eqnarray*}
We call a family $r^{\mathrm{ad}}_\e(\kappa), \kappa\in\Sigma$, such that
$u_\e^{\mathrm{ad}}\asymp u_\e(\kappa,r^{\mathrm{ad}}_\e(\kappa))$, the family of
\textit{adaptive separation rates}.

The relation $\g_\e(\Theta_{\e}(\Sigma))\to0$ is possible
if $u_\e(\Sigma)\to\infty$. It was shown in the supplementary material
\cite{ISS12}, Theorem
7.1, that
this relation suffices for the construction of minimax adaptive
consistent families of tests for mildly ill-posed inverse problems
with the class of analytic functions. This implication, however, does
not hold in the remaining ill-posed inverse problems under
consideration. In these cases, adaptive distinguishability
conditions and adaptive separation rates are sought, and they are
the goal of the subsequent sections. In contrast to the above mentioned theorem,
there is price to pay for the adaptation. We show that
$u_\e^{\mathrm{ad}}=\sqrt{\log\log\e^{-1}}$ for the mildly ill-posed inverse
problems with the Sobolev class of functions and
$u_\e^{\mathrm{ad}}={\log\log\e^{-1}}$ for other problems under consideration
(except the case mildly ill-posed inverse problems with the class of
analytic functions). These yield a loss in the separation rates
in terms of an extra $\sqrt[4]{\log\log\e^{-1}}$ factor for the mildly
ill-posed inverse problems with the Sobolev class, and in terms of an extra
$\sqrt{\log\log\e^{-1}}$ factor for severely problems with analytic
classes of functions. (A similar loss in the separation rates for a well-posed
signal detection problem was first observed in~\cite{Spok}.)

As we shall show below, the derived families of tests are of simple structure.
In particular, for the mildly ill-posed inverse problems with the
Sobolev class of
functions, these are of
the form
\[
\psi_\e^{\mathrm{ad}}=\one_{\{\sup_k{t_{\e,m_k}>H_k}\}},\qquad
m_k=2^k,\qquad H_k=\sqrt{C\log(k)},\qquad k \geq L,
\]
where $C>2$, for an integer-valued family $L=L_{\e}$, $L_{\e}\to
\infty
$, and
%
%e5.1 #&#
%
\begin{equation}\label{test2}
t_{\e, m}=\frac{1}{\sqrt{2 m}}\sum_{k=1}^{ m}
\bigl((y_k/\e)^2-1\bigr)
\end{equation}
are centered and normalized version of $\chi^2$-statistics that
correspond to the first $m$ observations.
%$(y_1,y_2,\ldots,y_m)$ in \nref{1.03}).
%The losses in the separation rates $r_{\e}^*$
%are in terms of a $\sqrt{\log\log(\e^{-1})}$ factor.

For the severely ill-posed inverse problems with the Sobolev class of
functions or the
class of analytic functions,
the derived families of tests are of the form
\begin{eqnarray*}
\psi_\e^{\mathrm{ad}}&=&\one_{\{\sup_k{|y_k|}>\e H_k\}},\qquad H_k=
\sqrt{2\log(k)},\qquad k < L,\\
H_k&=&\sqrt{C\log(k)},\qquad k \geq L,
\end{eqnarray*}
where $C>2$, for an integer-valued family $L=L_{\e}$, $L_{\e}\to
\infty$.
%and the losses in the separation rates $r_{\e}^*$ are in terms of a
%$\log\log(\e^{-1})$ factor.

Finally, for the severely ill-posed inverse problems with the
generalized analytic class of
functions, the derived tests are of the form
\[
\psi_\e^{\mathrm{ad}}=\one_{\{\sup_k{|y_k|}>\e T_{\e,k}\}},\qquad T_{\e,k}=
\max\bigl(T_\e,\sqrt{2\bigl(\log(k)+\log\log(k)\bigr)} \bigr)
\]
for a family $T_{\e} \to\infty$.
%, and the losses in the values of
%$u_{\e}$ are in terms of a $\log\log(\e^{-1})$ factor.

%In what follows, $\bar A$ denotes the complement of a set $A$.

%s5.1.1 #&#
\subsubsection{Mildly ill-posed inverse problems with the Sobolev class
of functions}
\label{Ad1}

Consider first the ``standard'' case $q=2$. Let $a_k = k^\alpha$ and
$\sigma_k
= k^\beta$, $k\in\N$, \mbox{$\alpha>0$}, $\beta>0$. Set $\kappa=(\a,\b
)$, and
let $\Sigma$ be a compact subset of $\R_+^2$. We show that, under a
weak assumption on the set $\Sigma$,
$
u_\e^{\mathrm{ad}}=\sqrt{\log\log(\e^{-1})}.
$
This corresponds to the adaptive separation rates
%
%e5.2 #&#
%
\begin{equation}
\label{adsepratF1} r^{\mathrm{ad}}_\e(\kappa)=\bigl(\e\sqrt[4]{
\log\log\bigl(\e^{-1}\bigr)}\bigr)^{4\a/(4\a+4\b+1)}.
\end{equation}

The rate optimal adaptive family of tests is of the following
structure. Take a collection $m_k=2^k$, $k\in\N$, $k\geq L=L_{\e}$,
for an integer-valued family $L_{\e} \to\infty$,
$L_{\e}=o(\log(\e^{-1}))$, and a family of test statistics $t_{\e,
m_k}$ of the form~(\ref{test2}). Consider the thresholds and tests
%
%e5.3 #&#
%
\begin{equation}
\label{adcas1F}\qquad H_k=\sqrt{C\log(k)},\qquad \CY_\e=\{y\dvtx
t_{\e, m_k} > H_k, \forall k\ge L_\e\},\qquad
\psi_\e=\one_{{\CY}_\e},
\end{equation}
where $C>2$. Denote also
%
%e5.4 #&#
%
\begin{equation}
\label{phsetF1} \phi(\kappa)=\frac{4}{4\a+4\b+1},\qquad \phi(\Sigma)=\bigl\{
\phi(\kappa)\dvtx \kappa\in\Sigma\bigr\}\subset(0,\infty).
\end{equation}

%Then, the following statement is true.

%th5.1 #&#
%
\begin{theorem}
\label{FanAdC1} %Consider the GSM \nref{ffgsm} and the hypothesis
%testing problem \nref{1.4}-\nref{ell} for $$.
Let $q=2, a_k =
k^\alpha$ and $\sigma_k = k^\beta$, $k\in\N$, $\alpha>0, \beta>0$.

\begin{longlist}[(a)]
\item[(a)]
(Lower bounds.) Let the set $\phi(\Sigma)$ given by~(\ref{phsetF1})
contain an interval
$[a,b],
0<a<b<\infty$. Then, there exists constant $d>0$ such that if $\limsup
_{\kappa\in\Sigma}u_\e(\kappa)/\sqrt{\log\log(\e^{-1})}\le d$, then
$\g_\e(\Theta_\e(\Sigma))\to1$.

\item[(b)] (Upper bounds.) For the family of tests $\psi_{\e}$ given by (\ref
{adcas1F}),
$\a(\psi_\e)=o(1)$, and there exists
constant $D=D(\Sigma)>0$ such that if
$\liminf_{\kappa\in\Sigma}u_\e(\kappa)/\break\sqrt{\log\log(\e^{-1})}>D$,
then $\b_\e(\psi_\e,\Theta_\e(\Sigma))=o(1)$.

\item[(c)]
(Adaptive separation rates.) The adaptive separation rates $r_{\e
}^{\mathrm{ad}}(\kappa)$,\break
$\kappa\in\Sigma$, are
given by~(\ref{adsepratF1}).
\end{longlist}
\end{theorem}

The proof is given in the supplementary material~\cite{ISS12}, Section 11.11.\vadjust{\goodbreak}

%re5.1 #&#
%
\begin{remark}
%In view of Remark~\ref{FFF1}, the rate optimality results obtained in
%Theorem
%asymptotic results and
%the adaptive separation rates \nref{adsepratF1} are true not only for
%the ill-posed
%problems ($\beta>0$) under consideration but also for a range of
%corresponding
%well-posed problems ($\beta\in(-1/4, 0]$). Furthermore, in view of
%Remark
%results and the
%same adaptive separation rates hold for the sequences $a_k \sim k^{
%$\sigma\sim k^{\beta}$, $k\in\N$, $\alpha>0,
%$k\in\N$, $\alpha>0,
%
%%% New revised remark %%%%
In view of Remark~\ref{RemequivSS}, similar rate optimality results
and the
same adaptive separation rates hold for the sequences $a_k \sim
k^{\alpha}$ and
$\sigma\sim k^{\beta}$, $k\in\N$, $\alpha>0,
\beta>0$, and the sequences $a_k \asymp k^{\alpha}$ and $\sigma
\asymp
k^{\beta}$,
$k\in\N$, $\alpha>0,
\beta>0$.
\end{remark}

As in the case of rate and sharp asymptotics, the ``sparse'' case $q
\in(0,2)$ is not
directly linked will Theorem~\ref{testing}. Rate-optimal adaptive tests
in this case, however, can also be constructed, based on the family of
tests considered in Section 7.4.1 of~\cite{IS02}, Chapter 7. Their
construction is omitted.
% and, hence, adaptive separation rates for the ``sparse'' case $q \in
%(0,2)$
% will be considered separately in Section~\ref{Ad1q}. (In fact,
%adaptivity results for the (more general) Besov class of functions
%will be presented therein.)

%s5.1.2 #&#
\subsubsection{Severely ill-posed inverse problems with the class of
analytic functions}
\label{Ad2} Let $a_k = \exp(\alpha k)$ and $\sigma_k = \exp(\beta
k)$, $k\in\N$, $\alpha>0, \beta>0$. Set $\kappa=(\a,\b, q)$,
and let
$\Sigma$ be a compact subset of $\R_+^2\times(0,2]$. We show that,
under a weak assumption on the set $\Sigma$,
$
u_\e^{\mathrm{ad}}=\log\log(\e^{-1}).
$
This corresponds to the adaptive separation rates
%
%e5.5 #&#
%
\begin{equation}
\label{adsepratF2} r^{\mathrm{ad}}_\e(\kappa)=\bigl(\e\sqrt{\log
\log\bigl(\e^{-1}\bigr)}\bigr)^{\a/(\a+\b)}.
\end{equation}

The rate optimal adaptive family of tests is of the following
structure. Take an
integer-valued family $L=L_{\e}$, $L_{\e} \to\infty$, $L_{\e
}=o(\log
\log(\e^{-1}))$.
Consider the families of thresholds and tests
%
%e5.6 #&#
%
\begin{equation}
\label{adcas2F}  H_k=\cases{\sqrt{2\log(L)},&\quad$k<L$,
\vspace*{2pt}\cr
\sqrt{C
\log(k)},&\quad$k\ge L$,}\qquad \CY_\e=\bigl\{y\dvtx |y_k| > \e
H_k\bigr\}, \quad\psi_\e=\one_{{
\CY}_\e},\hspace*{-35pt}
\end{equation}
where $C>2$. Denote also
%
%e5.7 #&#
%
\begin{equation}
\label{phsetF2} \phi(\kappa)=\frac{1}{2(\a+\b)},\qquad \phi(\Sigma)=\bigl\{
\phi(\kappa)\dvtx \kappa\in\Sigma\bigr\}\subset(0,\infty).
\end{equation}

%th5.2 #&#
%
\begin{theorem}
\label{FanAdC2} %Consider the GSM \nref{ffgsm} and the hypothesis
%testing problem \nref{1.4}-\nref{ell}.
Let $a_k = \exp(\alpha k), \sigma_k = \exp(\beta k), k\in\N,
\alpha>0, \beta>0$.

\begin{longlist}[(a)]
\item[(a)]
(Lower bounds.) Let the set $\phi(\Sigma)$ given by~(\ref{phsetF2})
contains an interval
$[a,b],
0<a<b<\infty$. Then, there exists constant $d>0$ such that if $\limsup
_{\kappa\in\Sigma}u_\e(\kappa)/\log\log(\e^{-1}) \leq d$, then
$\g_\e(\Theta_\e(\Sigma))\to1$.

\item[(b)] (Upper bounds.) For the family of tests $\psi_{\e}$ given by
(\ref{adcas2F}), $\a(\psi_\e)=o(1)$ and there exists constant
$D=D(\Sigma)>0$ such that if
$\liminf_{\kappa\in\Sigma}u_\e(\kappa)/\break\log\log(\e^{-1})>D$, then
$\b_\e(\Theta_\e(\Sigma),\psi_\e)=o(1)$.

\item[(c)]
(Adaptive separation rates.) The adaptive separation rates $r_{\e
}^{\mathrm{ad}}(\kappa)$,\break $\kappa\in\Sigma$, are
given by~(\ref{adsepratF2}).
\end{longlist}
\end{theorem}

The proof is given in the supplementary material~\cite{ISS12}, Section 11.12.

%re5.2 #&#
%
\begin{remark}
In view of Remark~\ref{RemequivSevAn}, similar rate optimality
results and the adaptive
separation rates~(\ref{adsepratF2}) hold for the sequences $a_k \asymp
\exp(\alpha k)$ and
$\sigma\asymp\exp(\beta k)$, $k\in\N$, $\alpha>0,
\beta>0$.
\end{remark}

%s5.1.3 #&#
\subsubsection{Severely ill-posed inverse problems with the Sobolev
class of functions}
\label{Ad3}

Let $a_k = k^\alpha$ and $\sigma_k = \exp(\beta k)$, $k\in\N$,
$\alpha>0, \beta>0$. Set $\kappa=(\a,\b, q)$, and let
$\Sigma$ be a compact subset of $\R_+^2\times(0,2]$. We show that, under
a weak assumption on the set $\Sigma$,
$
u_\e^{\mathrm{ad}}=\log\log(\e^{-1}).
$
This corresponds to the adaptive separation rates
%
%e5.8 #&#
%
\begin{eqnarray}
\label{adsepratF3} r_{\e}^{\mathrm{ad}}(\kappa)&=& \biggl(
\frac{2\beta}{2\log(\e^{-1})-2\alpha\log\log(\e^{-1})-\log\log
\log(\e^{-1})} \biggr)^{\alpha} \nonumber\\[-9pt]\\[-9pt]
&\sim&\biggl(\frac{\b}{\log(\e^{-1})}
\biggr)^{\alpha}.\nonumber
\end{eqnarray}

The rate optimal adaptive family of tests is of the following structure.
Take an integer-valued family $L=L_{\e}$, $L_{\e} \to\infty$,
$L_{\e
}=o(\log\log(\e^{-1}))$,
and consider the families of thresholds and tests given
by~(\ref{adcas2F}).\vspace*{-1pt}

%th5.3 #&#
%
\begin{theorem}
\label{FanAdC3} %Consider the GSM \nref{ffgsm}
%and the hypothesis testing problem \nref{1.4}-\nref{ell}.
Let $a_k = k^{\alpha}$ and
$\sigma_k = \exp(\beta k)$, $k\in\N$, $\alpha>0, \beta>0$.

\begin{longlist}[(a)]
\item[(a)]
(Lower bounds.) Let the set $\Sigma$ contains an interval of
$(\a,\b)\dvtx \b\in[1/2b$, $1/2a]$, $0<a<b<\infty$, and a fixed $\a>0$.
Then there exists constant $d>0$ such that if $\limsup_{\kappa\in\Sigma
}u_\e(\kappa)/\log\log(\e^{-1}) \leq d$, then
$\g_\e(\Theta_\e(\Sigma))\to1$.

\item[(b)] (Upper bounds.) For the family of tests $\psi_{\e}$ given by (\ref
{adcas2F}),
$\a(\psi_\e)=o(1)$ and there exists
constant $D=D(\Sigma)>0$ such that if
$\liminf_{\kappa\in\Sigma}u_\e(\kappa)/\break \log\log(\e^{-1})>D$, then
$\b_\e(\psi_\e,\Theta_\e(\Sigma))=o(1)$.

\item[(c)]
(Adaptive separation rates.) The adaptive separation rates $r_{\e
}^{\mathrm{ad}}(\kappa)$,\break $\kappa\in\Sigma$, are
given by~(\ref{adsepratF3}).\vspace*{-1pt}
\end{longlist}
\end{theorem}

The proof is given in the supplementary material~\cite{ISS12}, Section 11.13.\vspace*{-1pt}

%re5.3 #&#
%
\begin{remark}
\label{adNewFF2S1}
It is worth mentioning that a stronger result is possible in this case.
In view of (11.21) in~\cite{ISS12}, relation
(\ref{adsepratF3})
determines \textit{sharp adaptive separation rates} $r_{\e}^{\mathrm{ad}}(\kappa)$,
$\kappa\in\Sigma$, in the
following sense:

\begin{longlist}[(a)]
\item[(a)]
if $\liminf(r_\varepsilon(\kappa) / r_{\e}^{\mathrm{ad}}(\kappa
) )
>1$, then $u_\varepsilon\to\infty$, that is,
$\gamma_\varepsilon(r_\varepsilon)\to0$;

\item[(b)] if $\limsup(r_\varepsilon(\kappa) / r_{\e}^{\mathrm{ad}}(\kappa
) )<1$,
then $u_\varepsilon\to0$, that is,
$\gamma_\varepsilon(r_\varepsilon)\to1$, and the minimax testing
is impossible.\vspace*{-1pt}
\end{longlist}
\end{remark}

%re5.4 #&#
%
\begin{remark}
In view of Remark~\ref{equivR3FF2}, similar rate optimality
results and the adaptive separation rates~(\ref{adsepratF3}) (as well
as the sharp adaptive separation rates mention in Remark~\ref{adNewFF2S1})
hold for the sequences $a_k \sim k^{\alpha}$ and $\sigma\asymp
\exp(\beta k)$, $k\in\N$, $\alpha>0, \beta>0$.\vspace*{-1pt}
\end{remark}

%s5.1.4 #&#
\subsubsection{Extremely ill-posed inverse problems with
the class of generalized analytic functions}\label{Ad5}

We consider the case $q=2$ only. By the results of Section~\ref{R5},
in order to obtain distinguishability conditions, we can replace
$u_\e(\kappa)$, $\kappa\in\Sigma$, by
$u_\e^{\mathrm{lin}}(\kappa)=u_{\e}^{\mathrm{lin}}(\kappa,r_\e(\kappa))$,
determined\vspace*{1pt}
by~(\ref{elinFFF}), by $a_k=a_k(\kappa)$ and
$\sigma_k=\sigma_k(\kappa)$, $k\in\N$. Set $
u_\e^{\mathrm{lin}}(\Sigma)=\inf_{\kappa\in\Sigma}u_\e^{\mathrm{lin}}(\kappa). $ We
are interested in how large $u_\e^{\mathrm{lin}}(\Sigma)$ should be in order
to $\g_\e(\Theta_{\e}(\Sigma))\to0$.\vadjust{\goodbreak}

Assume below the uniform version of~(\ref{assump}): let
$\sigma_k(\kappa)$ and $a_k(\kappa)$, $k\in\N$, be increasing
sequences such that, for all $\kappa\in\Sigma$ and some constants
$0<b<B$,
%
%e5.9 #&#
%
\begin{equation}
\label{assuni0} b\le a_1(\kappa)\le B,\qquad b\le
\sigma_1(\kappa)\le B
\end{equation}
and, for some increasing sequence $\tau_k>1, \tau_k\to\infty$ and
some $c_0>1$ for all $\kappa\in\Sigma$ and $k\in\N$,
%
%e5.10 #&#
%
\begin{equation}
\label{assuni} \sigma_{k+1}(\kappa)/\sigma_k(
\kappa)\ge\tau_k,\qquad a_{k+1}(\kappa)/a_k(\kappa)\ge
c_0.
\end{equation}

Similar to $u_\e^{\mathrm{ad}}$, one can consider a family
$u_{\e,\mathrm{ad}}^{\mathrm{lin}}$ which characterizes adaptive distinguishability.
We show that, under some assumption on the set $\Sigma$, one has
$u_{\e,\mathrm{ad}}^{\mathrm{lin}}=\log\log(\e^{-1})$.

For $\kappa\in\Sigma$ and for $A>0$ large enough, let an integer
$m=m(A,\kappa)$ be defined by the relations
%
%e5.11 #&#
%
\begin{equation}
\label{mA} a_{m-1}(\kappa)\sigma_{m-1}(\kappa)\le A<
a_{m}(\kappa)\sigma_{m}(\kappa).
\end{equation}
Under~(\ref{assuni0}) and~(\ref{assuni}), one has
$
a_{m-1}\sigma_{m-1}\ge b^2c_0^{m-2}\prod_{k=1}^{m-2}\tau_k,
$
which yields
%
%e5.12 #&#
%
\begin{equation}
\label{MA1} \sup_{\kappa\in\Sigma}m(A,\kappa)=o\bigl(\log(A)\bigr)
\qquad\mbox
{as } A\to
\infty.
\end{equation}
Set
$
\CM(A,\Sigma)=\{m(A,\kappa)\in\N\dvtx \kappa\in\Sigma\},
M(A,\Sigma)=\#(\CM(A,\Sigma)).
$
Since $M(A$, $\Sigma)\le\max_{m\in\CM(A,\Sigma)}m$, one has, by
(\ref{MA1}),
$M(A,\Sigma)=o(\log(A))$ as $A\to\infty$.
Let $m=m(A,\kappa)$ be defined by~(\ref{mA}) and set $L(A,\Sigma):=
\sup_{\kappa\in\Sigma}\log(m(A,\kappa))$. By~(\ref{MA1}) we
have, as
$A\to\infty$,
%
%e5.13 #&#
%
\begin{equation}
\label{mmin} %m(A,\Sigma)\eq\inf_{\kappa\in\Sigma}m(A,\kappa)\to\infty,
\limsup L(A,\Sigma)/\log\log(A)\le1.
\end{equation}
%
%Set %$m_\e(\Sigma)=m(A_\e,\Sigma)$ for
%$A_\e=(\e\sqrt{\log\log(\e^{-1})})^{-1}$.

For the lower bounds we suppose one can find %a subset $\tilde\Sigma
quantities $b>0, C\ge1$ such that
%
%e5.14 #&#
%
\begin{eqnarray}
\label{mAA}
\liminf_{A\to\infty}\log\bigl(M(A,\Sigma)\bigr)/\log\log(A)&=&b,\nonumber\\[-8pt]\\[-8pt]
\sup_{\kappa\in\Sigma} u_\e^{\mathrm{lin}}\bigl(\kappa,r_\e(
\kappa)\bigr)&\le& C u_\e^{\mathrm{lin}}(\Sigma).\nonumber
\end{eqnarray}
[The first relation in~(\ref{mAA}) is fulfilled for the example
mentioned in
Remark~\ref{remexample}, at least if the set
$\Sigma=\{(\a,\tau,\b,\gamma)\}$ contains an interior point.]

The rate optimal adaptive family of tests is of the following
structure. Take a family $T_\e\to\infty$ such that
$T_\e=o(\sqrt{\log\log(\e^{-1})})$, and take a family of sequences
$T_{\e,k}$ of the form $
T_{\e,k}=\max(T_\e,\sqrt{2(\log(k)+\log\log(k))} )$.

Consider the families of events and of tests
%
%e5.15 #&#
%
\begin{equation}
\label{adcas5F} \CY_{\e}=\bigl\{y\dvtx |y_k| > \e
T_{\e,k}, \forall k\in\N\bigr\},\qquad \psi_\e=\one_{{\CY}_\e}.
\end{equation}

%Then, the following statement is true.

%th5.4 #&#
%
\begin{theorem}
\label{FanAdC5} Consider the GSM~(\ref{ffgsm})
and the hypothesis testing problem~(\ref{14}) and~(\ref{ell}) for $q=2$.
Let $\{a_k\}_{k\in\N}$ and
$\{\sigma_k\}_{k\in\N}$ be increasing sequences satisfying
(\ref{assuni0}) and~(\ref{assuni}).

\begin{longlist}[(a)]
\item[(a)]
(Lower bounds.) Assume~(\ref{mAA}). Then there exists a constant
$d>0$ such that if $\lim\sup u_\e^{\mathrm{lin}}(\Sigma)/ \log\log(\e^{-1})
\leq d$, then $\g_\e(\Theta_\e(\Sigma))\to1$.

\item[(b)]
(Upper bounds.) Assume~(\ref{mmin}). For the family of tests
$\psi_{\e}$ given by~(\ref{adcas5F}), there exists a constant $D>0$
such that if $u_\e^{\mathrm{lin}}(\Sigma)>D\log\log(\e^{-1})$, then
$\g_\e(\Theta_\e(\Sigma), \psi_\e)=o(1)$.

\item[(c)]
(Adaptive separation rates.) The adaptive separation rates
$r^{\mathrm{ad}}_\e(\kappa),\break \kappa\in\Sigma$, are determined by the relation
$u_{\e, \mathrm{ad}}^{\mathrm{lin}}\asymp u_\e^{\mathrm{lin}}(\kappa,r^{\mathrm{ad}}_\e(\kappa))$.
\end{longlist}
\end{theorem}

The proof is given in the supplementary material~\cite{ISS12}, Section 11.14.

\section*{Acknowledgments}

We would like to thank the Editor, the Associate Editor and two
anonymous referees for useful comments and suggestions.

\begin{supplement}%[id=suppA]
\stitle{Detailed proofs and other material}
\slink[doi]{10.1214/12-AOS1011SUPP} %[doi,text={...}] - jei reikia
%suskaldyti doi
\sdatatype{.pdf}
\sfilename{aos1011\_supp.pdf}
\sdescription{In this supplement, we present relevant material and the
detailed proofs of the previous sections.}
\end{supplement}

% imsref loaded by lrinkeviciute, 2012-07-04 11:21:09
% imsref loaded by lrinkeviciute, 2012-07-04 12:14:11
%

\printaddresses

\end{document}